\documentclass[11pt,a4paper]{article}

\usepackage[english]{babel}
\usepackage{amssymb}
\usepackage[latin1]{inputenc}
\usepackage{graphicx}
\usepackage{amsthm}
\usepackage{color}
\usepackage{stmaryrd}
\usepackage{fullpage}

\title{The strong inviscid limit of the isentropic compressible Navier-Stokes equations with Navier boundary conditions}
\date{October 2014}
\author{\textsc{Matthew Paddick} \\ LJLL, Universit\'e Pierre \& Marie Curie, Paris, France \\ \texttt{paddick@ljll.math.upmc.fr}}

\begin{document}

\renewcommand{\labelitemi}{$\bullet$}
\newtheorem{theo}{Theorem}[section]
\newtheorem{lemma}[theo]{Lemma}
\newtheorem*{nota}{Notation}
\newtheorem{assu}[theo]{Assumption}
\newtheorem{propo}[theo]{Proposition}
\newtheorem{coro}[theo]{Corollary}
\newtheorem{esti}[theo]{Proposition}
\newcommand{\preu}[1]{\textit{\underline{#1}}}
\newcommand{\thref}[1]{Theorem \ref{#1}}
\newcommand{\lemref}[1]{Lemma \ref{#1}}
\newcommand{\propref}[1]{Proposition \ref{#1}}
\newcommand{\cororef}[1]{Corollary \ref{#1}}
\newcommand{\estiref}[1]{Proposition \ref{#1}}
\newcommand{\assuref}[1]{Assumption \ref{#1}}
\newcommand{\tworef}[3]{#3 \ref{#1} and \ref{#2}}

\newcommand{\til}[1]{\tilde{#1}}
\newcommand{\rplus}{\mathbb{R}^{+}}
\newcommand{\rtwo}{\mathbb{R}^{2}}
\renewcommand{\div}{\mathrm{div}~}
\newcommand{\divt}{\mathrm{div}\tang}
\newcommand{\rot}{\overrightarrow{\mathrm{rot}~}}
\newcommand{\eps}{\varepsilon}
\newcommand{\norme}[2]{\left\| #2 \right\| _{#1}}
\newcommand{\trinorme}[2]{\interleave #2 \interleave _{#1}}
\newcommand{\tang}{_\tau}
\newcommand{\derp}[1]{\partial_{#1}}
\newcommand{\diag}{\mathrm{diag}}
\newcommand{\dOmega}{\partial\Omega}
\newcommand{\Lip}{\mathrm{Lip}}
\newcommand{\rhs}[1]{\mathbf{RHS}(\ref{#1})}
\newcommand{\mathand}{~~\mathrm{and}~~}
\newcommand{\mathif}{\mathrm{if}~}
\newcommand{\sdel}[1]{#1_{\delta}}

\maketitle

\begin{abstract}
We obtain existence and conormal Sobolev regularity of strong solutions to the 3D compressible isentropic Navier-Stokes system on the half-space with a Navier boundary condition, over a time that is uniform with respect to the viscosity parameters when these are small. 
These solutions then converge globally and strongly in $L^2$ towards the solution of the compressible isentropic Euler system when the viscosity parameters go to zero.
\newline


\end{abstract}

\section{Introduction}

\subsection{The isentropic compressible Navier-Stokes system and the inviscid limit problem}

We consider the motion of a compressible fluid in the half-space $\Omega=\rtwo\times\rplus$. The density of the fluid is a scalar function $\rho(t,x)$ with $t$ being the time variable and $x \in\Omega$, its velocity is $u(t):\Omega \rightarrow \mathbb{R}^{3}$: these will be the unknowns of the system. 
The temperature of the fluid, $\theta$, is constant throughout the paper, and the pressure $P$ will follow a barotropic law: $P=P(\rho)=\frac{k}{\gamma}\rho^{\gamma}$, with $k>0$ and the adiabatic constant $\gamma> 1$. 
The motion of the fluid is governed by the isentropic compressible Navier-Stokes system, which consists of two conservation laws:
\begin{itemize}
\item conservation of mass
\begin{equation} \derp{t}\rho + \div(\rho u) = 0, \label{NSdiv} \end{equation}
\item and conservation of momentum
\begin{equation} \derp{t}(\rho u) + \div(\rho u\otimes u) = \div \Sigma + \rho F, \label{NS} \end{equation}
\end{itemize}
with $F(t,x)$ a force term. To avoid technical complications with compatibility conditions, we consider that the force is smooth and that fluid is at rest for negative times: $F(t,x)=0$ and $(\rho,u)(t,x)=(1,0)$ for $t<0$.

In these equations, $\Sigma$ is the internal stress tensor,
$$ \Sigma = 2\eps\mu Su + (\eps\lambda_0 \div u - P(\rho)) I_{3}, $$
in which $\eps>0$ will be arbitrarily small, $Su=\frac{1}{2}(\nabla u+\nabla u^{T})$ is the symmetric gradient, and $I_{3}$ is the $3\times 3$ identity matrix. The parameters $\mu>0$ (dynamic viscosity) and $\lambda:=\lambda_0 + \mu >0$ (bulk viscosity) will be given regular functions of $\rho$. In equation (\ref{NS}), we will write
$$ \div \Sigma = \eps(\mu \Delta u +\lambda \nabla\div u) - \nabla (P(\rho)) + \eps\sigma(\nabla \rho,\nabla u). $$
The final term $\sigma$ is non-zero only when $\lambda'(\rho)$ and $\mu'(\rho)$ are non-zero ($\lambda$, $\mu$ not constant).
\newline

We shall not go into detailed historics about the isentropic system; a survey of existence and regularity results up to 1998 was put together by B.~Desjardins and C-K.~Lin, \cite{DL}. 
We shall quickly cite the emblematic result on the subject: the global existence theorem for weak solutions, as proved by P-L.~Lions in the 1990s (\cite{Lpl} or \cite{Lplbook2}), a result that has since been improved, for instance, by E.~Feireisl, A.~Novotn\'y and H.~Petzeltov\'a, \cite{FNP}, 
and extended to some cases with density-dependent (variable, therefore) viscosity coefficients by D.~Bresch, B.~Desjardins and D.~G\'erard-Varet \cite{BDG}. 
In the latter paper, a nonlinear drag force, written as $F=F(u)=r_0 |u| u$, is considered. 
We also refer to the local strong existence theory, initiated in the 1960s and 70s by J.~Nash and then V.A.~Solonnikov \cite{Nj, Sva} and improved on by R.~Danchin (for instance, see survey \cite{Drsurv}). In this framework, blowup can occur \cite{Xz}.
\newline

In this article, we will aim to obtain local-in-time existence of strong solutions to the isentropic Navier-Stokes system, with a lower bound for the time of existence that does not depend on the viscosity parameters when these are small, and work on the inviscid limit problem. 
As we are on the half-space, we will have some boundary conditions.

\begin{nota} Throughout the paper, we will use the notation $x=(y,z)$, with $y\in\rtwo$ and $z\in\rplus$; the boundary is therefore the set $\{x=(y,z)\in\mathbb{R}^{3}~|~z=0\}$.

Moreover, for a vector field $v(x)$, the tangential part to the boundary is, for $x$ on the boundary, $v\tang(x)=v(x)+(v(x)\cdot\vec{n}(x))\vec{n}(x)$, where $\vec{n}(x)$ is the outer normal vector to $\dOmega$ at point $x$. As $\vec{n}(x)\equiv(0,-1)$, we extend the notation to all $\Omega$: $v\tang(x)=(v_{1}(x),v_{2}(x))$. \end{nota}

The boundary conditions on $u$ are the standard non-penetration condition on the boundary,
\begin{equation} u\cdot\vec{n}|_{z=0}=u_{3}|_{z=0}=0 , \label{NPB} \end{equation}
and the Navier (slip) boundary condition $\left[\left(\frac{1}{\eps}\Sigma\vec{n}+ au\right)\tang\right]|_{z=0} = 0$ in which $a>0$. In our case, with a flat boundary, the Navier condition can be rewritten as
\begin{equation} [\mu(\rho) \derp{z}u\tang]|_{z=0} = 2a u\tang |_{z=0} \label{NBC} \end{equation}
for $t>0$ and $y\in\rtwo$. As opposed to the Dirichlet or no-slip condition, which, in our setting, would be $u|_{z=0}=0$, the Navier condition, proposed by H.~Navier himself in the XIX\textsuperscript{th} century \cite{Nh}, allows the fluid to slip along the boundary, and this occurs wherever interaction at the boundary is non negligeable. 
For instance, the slip phenomenon can be observed on the contact line of two immiscible flows \cite{QWS}, and in capillary blood vessels, which are the microscopic, tissue-irrigating vessels where molecular exchanges with the neighbouring cells take place \cite{PRD}. 
It also appears when homogenising rough and porous boundaries (\cite{JM} and \cite{GVM}), and can be derived mathematically from a Boltzmann microscopic model with a Maxwell reflection boundary condition \cite{MSR}. 
To be physically pertinent, the slip coefficient $a$ should be chosen positive, but our results do not technically require $a$ to have a specific sign, so we take $a\in \mathbb{R}$.

We also impose the limit condition
\begin{equation} U(t,x):=(\rho(t,x)-1,u(t,x))\stackrel{|x|\rightarrow+\infty}{\rightarrow} 0 \label{limcond} \end{equation}
so that $U(t) \in L^{2}(\Omega)$.
\newline

Formally, taking $\eps=0$ leads to the compressible Euler equations
\begin{equation} \left\{ \begin{array}{rcl} \derp{t}\rho + \div(\rho u) & = & 0 \\
\rho \derp{t}u + \rho u\cdot\nabla u + \nabla P(\rho) & = & \rho F . \end{array} \right. \label{Euler} \end{equation}
The Euler equation is of order one, so it only requires one boundary condition equation, which is (\ref{NPB}). 
This leads to the appearance of boundary layers: if the solutions of low-viscosity Navier-Stokes equations are expected to behave like a solution of the Euler equation far away enough from the boundary, solutions of Navier-Stokes are still required to satisfy a second boundary condition, whereas the reference solution of the Euler equation is not. 
A typical boundary layer expansion for solutions to Navier-Stokes will read
\begin{equation} u^\eps(t,y,z) = u^E(t,y,z) + V \left(t,y,\frac{z}{\sqrt{\eps}}\right) , \label{blayer} \end{equation}
with $u^E$ solving the Euler equation, and $V$, acting on a shorter scale, picking up the boundary condition.

The inviscid limit problem is a major challenge for mathematicians, whether one considers compressible or incompressible fluids. We remind the reader of the (non-forced) incompressible system:
$$ \left\{ \begin{array}{rcl} \div u & = & 0 \\
\derp{t}u + u\cdot\nabla u - \eps\nu\Delta u + \nabla q & = & 0 \\ u|_{t=0} & = & u_0 , \end{array} \right. $$
in which $\nu=\frac{\mu}{\rho}$, with $\rho$ constant, is the kinematic viscosity and $q$ is the kinematic pressure. 
Regarding the inviscid limit results on the incompressible system with Navier boundary conditions, the problem is solved in $L^2$ framework in 2D (see \cite{Bc}, \cite{CMR}, \cite{Kjp}), and convergence of weak solutions of the Navier-Stokes equation towards a strong solution of the Euler equation, when the limit initial condition is regular enough, has been obtained for a range of Navier slip coefficients of the form $a=a'\eps^{-\beta}$: 
starting with D.~Iftimie and G.~Planas \cite{IP} ($\beta=0$), \cite{WWX} ($\beta<1/2$) and \cite{Pm11} ($\beta<1$ for positive slip coefficients and $\beta\leq 1/2$ regardless of sign) have extended the range of numbers $\beta$ for which convergence occurs. 
C.~Bardos, F.~Golse and L.~Paillard recently obtained weak convergence results for Leray solutions of Navier-Stokes towards dissipative solutions of the Euler equation \cite{BGP}.

The solutions to the 3D incompressible Navier-Stokes equations with a Navier boundary condition that does not depend on $\eps$ have a better asymptotic expansion than (\ref{blayer}); D.~Iftimie and F.~Sueur showed in \cite{IS} that the boundary layer $V$ has a smaller amplitude:
\begin{equation} u^\eps(t,y,z)=u^E(t,y,z) + \sqrt\eps V(t,y,\eps^{-1/2}z) . \label{blayern} \end{equation}
The problems of local existence of strong solutions on a time interval that does not depend on $\eps$, and of the corresponding inviscid limit, showing behaviour in agreement with this ansatz, have been solved by N.~Masmoudi and F.~Rousset \cite{MR}. 
Their approach is based on energy estimates in conormal Sobolev spaces, and the same technique has allowed them to prove similar results for the corresponding free-boundary system \cite{MRfs}, and we will see that a similar approach for the isentropic system is valid. We refer the reader to \cite{XX} and \cite{BC2} for other studies in 3D.

In the incompressible case also, a variable viscosity coefficient, $\nu=\nu(q,|Su|^2)$, can be considered. This appears in elastohydrodynamics or the mechanics of granular materials for example. 
Existence of weak solutions for such equations on a bounded domain with the Navier boundary condition has been obtained by M.~Bul\'i\v{c}ek, J.~M\'alek and K.~Rajagopal \cite{BMR}.
\newline

On the 3D isentropic Navier-Stokes system with Navier boundary conditions that we are interested in, F.~Sueur recently showed the convergence of weak solutions to a strong solution of the Euler equation when the limit initial condition is smooth, and for slip coefficients that can depend on $\eps$, such as $a/\eps^\beta$ with $\beta<1$ \cite{Sf13}. 
In the case where the slip coefficient does not depend on $\eps$, D.~Hoff obtained global solutions with intermediate regularity (more regularity than weak solutions \textit{\`a-la-}Lions, but not classical solutions) in 2005 \cite{Hoff}, while Y-G.~Wang and M.~Williams justified a WKB expansion for strong solutions \cite{WW} in 2012. 
In particular, Wang and Williams show that solutions to the isentropic Navier-Stokes system behave similarly to their incompressible counterparts in the inviscid limit, in the sense that we have the asymptotic expansion (\ref{blayern}). 
As we aim to get existence of strong solutions of the Navier-Stokes equation through uniform \textit{a priori} estimates, we will make use of conormal derivatives that we introduce below.
\newline

To complete the references on the inviscid limit problem for compressible fluids, we cite F.~Huang, Y.~Wang and T.~Yang (ideal gas, \cite{HWY}), and Feireisl and Novotn\'y (weak-strong convergence in the Navier-Stokes-Fourier setting, \cite{FN13}), for advances on the full compressible Navier-Stokes equations, with an extra equation on the internal energy, temperature or entropy. 
We also refer to \cite{XY} for results on the linearised 2D system, \cite{Rf05} and \cite{GMWZ} for boundary layer analysis with characteristic and non-characteristic boundary conditions respectively, and \cite{MZ}, \cite{AB} for results on more general parabolic-hyperbolic systems.

\subsection{The conormal functional setting}

If we consider a boundary-layer expansion of the form (\ref{blayern}), which the solutions of the equation we will study satisfy, we see that we can expect uniform control of $u^\eps$, its derivatives, but not its second derivatives: a factor $\eps^{-1/2}$ results from the differentiation of the boundary layer $V$. 
Thus, the functional setting used in this paper will be that of conormal Sobolev spaces. Introduced in the mid-60s \cite{Hl}, these spaces have been used to work on hyperbolic systems with characteristic boundaries (see, for example, \cite{Rj}, \cite{Go}, \cite{Sp95}). 
Such spaces on a domain $\Omega$, which has a boundary, are constructed by differentiating functions following a finite set of generators of vector fields that are tangent to the boundary of $\Omega$. Namely, in the case of the half-space, we can choose
$$ Z_{1,2} = \derp{y_{1},y_{2}},~Z_{3} = \phi(z)\derp{z} , $$
with $\phi$ a smooth, positive, bounded function of $\rplus$ such that $\phi(0)=0$ and $\phi'(0)\neq 0$ - typically, consider $\phi(z)=\frac{z}{1+z}$. 

Conormal derivatives will allow us to get high-order uniform-in-$\eps$ estimates. Considering an expansion like (\ref{blayern}), if we look at the conormal derivatives of $\derp{z}u^\eps$, the boundary-layer term is written as
$$ Z^3 (\partial V(\eps^{-1/2}z)) = \eps^{-1/2}\phi(z)\partial^2 V(\eps^{-1/2}z) , $$
and this is of amplitude ${\cal O}(1)$ in a neighbourhood of the boundary of size $\sqrt\eps$ thanks to the factor $\phi(z)$. Thus the conormal setting is the only one in which we can expect uniform bounds on a large number of derivatives.

The conormal Sobolev space on $\Omega$, $W^{m,p}_{co}(\Omega)$, is then naturally defined as the set of functions $f(x) \in L^p(\Omega)$ such that the conormal derivatives of order at most $m$ of $f$ are also in $L^p(\Omega)$. 
\newline

As part of their estimation process in \cite{MR}, Masmoudi and Rousset used the incompressibility equation to express $\derp{z}u_3$ as a combination of conormal derivatives:
$$ \derp{z}u_3 = -\derp{y_1}u_1 -\derp{y_2}u_2 = -Z_1 u_1 - Z_2 u_2. $$
We will be able to use a similar trick, but with the equation of conservation of mass (\ref{NSdiv}), in which $\derp{t}\rho$ intervenes. Also, we will regularly use equation (\ref{NS}) to replace terms with two normal derivatives ($\derp{zz}u$), and there, $\derp{t}u$ is involved. 
For these reasons, we add $Z_0 = \derp{t}$ for functions that depend on $(t,x)$, and introduce the conormal Sobolev spaces on $[0,T]\times \Omega$ in the sense of O.Gu\`es \cite{Go}, for a set time $T$. 
For $\alpha\in\mathbb{N}^{4}$, we write $Z^{\alpha}=Z_{0}^{\alpha_{0}}Z_{1}^{\alpha_{1}}Z_{2}^{\alpha_{2}}Z_{3}^{\alpha_{3}}$, and $|\alpha|=\sum_{i=0}^{3} \alpha_i$: the conormal Sobolev space $W^{m,p}_{co}([0,T]\times \Omega)$ is the set of functions $f:
[0,T]\times \Omega\rightarrow\mathbb{R}^{d}$ such that $Z^{\alpha}f\in L^{p}([0,T]\times \Omega)$, for every $\alpha$ with $|\alpha|\leq m$. We will only use $p=+\infty$ and $p=2$, with the notation $H^m_{co}=W^{m,2}_{co}$. We therefore have
$$ H^m_{co}([0,T]\times\Omega) = \{ f(t,x)~|~ \forall ~0\leq k\leq m, ~ \derp{t}^k f \in L^2([0,T],H^{m-k}_{co}(\Omega))\} . $$
Compared to $W^{m,p}_{co}(\Omega)$, the notation is slightly abusive, in that $W^{m,p}_{co}([0,T]\times\Omega)$ is not a space whose conormal derivatives are tangent to the boundary ($Z_0=\derp{t}$ is not tangent to the boundary of $[0,T]$).

In our \textit{a priori} estimation process, we will be interested in the following space:
$$ X^m_T(\Omega) = \{ f(t,x) ~|~ \forall ~0\leq k\leq m,~ \derp{t}^k f\in L^\infty([0,T],H^{m-k}_{co}(\Omega)\}. $$
This is more restrictive than asking for $f\in H^{m}_{co}([0,T]\times\Omega)$. For a set $t\geq 0$, we introduce the semi-norms
$$ \norme{m}{f(t)}^2 = \sum_{|\alpha|\leq m} \norme{L^2(\Omega)}{Z^\alpha f(t)}^2 \mathand \norme{m,\infty}{f(t)} = \sum_{|\alpha|\leq m} \norme{\infty}{Z^\alpha f(t)}. $$
Note that these semi-norms coincide with the $H^m_{co}(\Omega)$ and $W^{m,\infty}_{co}(\Omega)$ norms if $f$ is stationary. Based on these semi-norms, we construct two norms on $X^m_T$ and $W^{m,\infty}_{co}([0,T]\times\Omega)$ respectively, which are essentially $L^\infty$-in-time norms,
$$ \trinorme{m,T}{f} := \sup_{t\in[0,T]} \norme{m}{f(t)} \mathand \trinorme{m,\infty,T}{f} := \sup_{t\in[0,T]} \norme{m,\infty}{f(t)} , $$
the latter of which coincides with the $W^{m,\infty}_{co}([0,T]\times\Omega)$ norm, and $L^2$-in-time norms
$$ \int_0^T \norme{m}{f(t)}^2~dt \mathand \int_0^T \norme{m,\infty}{f(t)}^2~dt , $$
the former of which is the natural norm for $H^m_{co}([0,T]\times\Omega)$. We will prefer not to abbreviate these last norms, as the forms we have given will make the \textit{a priori} estimation process clearer.

We add the following abbreviations: $\norme{\infty}{f(t)}:=\norme{0,\infty}{f(t)}$, $\trinorme{\infty,T}{f}:=\trinorme{0,\infty,T}{f}$, and we denote by $\trinorme{\Lip,T}{f}$ the standard Lipschitz norm on $[0,T]\times\Omega$.

\subsection{Results and proof strategy}

We introduce the notation $U(t,x)=(\rho-1,u)(t,x)$, and will consider the class of solutions satisfying the following property:
$$ {\cal E}_m(T,U) := \trinorme{m,T}{U}^2 + \trinorme{m-1,T}{\derp{z}u_{\tau}}^2+\int_0^T \norme{m-1}{\derp{z}u_{3}(s)}^2+\norme{m-1}{\derp{z}\rho(s)}^2~ds \hspace{50pt} $$
\begin{equation} \hspace{130pt} +\trinorme{1,\infty,T}{\derp{z}u_{\tau}}^2+\int_0^T \norme{1,\infty}{\derp{z}\rho(s)}^2+\norme{1,\infty}{\derp{t}\derp{z}\rho(s)}^2~ds <+\infty . \label{emtu} \end{equation}
In terms of our functional setting in the previous section, if we have $U\in X^m_T$ and $\derp{z} U\in X^{m-1}_T \cap W^{2,\infty}_{co}([0,T]\times\Omega)$ (with $L^\infty$ in time norms only), then ${\cal E}_m(T,U)$ is also finite. 

Note that in ${\cal E}_m(T,U)$ we only have control of $L^2$-in-time norms on the derivatives of $u_3$ and $\rho$. Simply using $|f(t,x)|^2 = |f(0,x)|^2 + \int_0^t \derp{t}(|f(s,x)|^2)~ds$, we get that
\begin{equation} \trinorme{1,\infty,T}{\derp{z}\rho}^2 \leq \norme{1,\infty}{\derp{z}\rho(0)}^2 + C\int_0^t \norme{1,\infty}{\derp{z}\rho(s)}^2+\norme{1,\infty}{\derp{t}\derp{z}\rho(s)}^2~ds , \label{l2timeinf} \end{equation}
thus the final integral in ${\cal E}_m(T,U)$ acts to control the $W^{1,\infty}_{co}$ norm of $\derp{z}\rho$. 
Requiring control of the $W^{1,\infty}_{co}$ norm of $\nabla U$ is typical of characteristic hyperbolic problems, see \cite{Go}.

The force must be smooth on $\mathbb{R}\times \Omega$, hence we introduce the notation:
$$ {\cal N}_m(T,F) = \sup_{t \in [-T,T]} \norme{m}{F(t)}^2 + \norme{m-1}{\nabla F(t)}^2 + \norme{2,\infty}{\nabla F(t)}^2 . $$

\begin{theo} \underline{Uniform existence of solutions to the Navier-Stokes system.}

\noindent Let $m\geq 7$, $F=(0,F)$ be such that ${\cal N}_m(t,F)<+\infty$ for any $t>0$, and $\mu$ and $\lambda$ be positive and bounded ${\cal C}^m$ functions of $\rho$. 
Then, for $\eps_0>0$, there exists $T^*>0$ such that, for every $0<\eps<\eps_0$, there is a unique $U^\eps$ satisfying ${\cal E}_m(T^*,\cdot)<+\infty$, solution to (\ref{NSdiv})-(\ref{NS})-(\ref{NPB})-(\ref{NBC})-(\ref{limcond}), the isentropic compressible Navier-Stokes system on $(0,T^*)\times\Omega$ with Navier boundary conditions. 
Moreover, there is no vacuum on this time interval: there exists $c>0$ such that $\rho(t,x)\geq c$ for $t\in [0,T^*]$ and $x\in\Omega$. \label{unifexist} \end{theo}

\begin{theo} \underline{Inviscid limit.}

Under the same conditions as above, the family $(U^\eps=(\rho^\eps-1,u^\eps))_{0<\eps<\eps_0}$ of solutions to the Navier-Stokes system converges in $L^2([0,T^*]\times\Omega)$ and $L^\infty([0,T^*]\times\Omega)$, towards $V=(\rho-1,u)$, the unique solution to the isentropic compressible Euler system, (\ref{Euler})-(\ref{NPB})-(\ref{limcond}), that satisfies ${\cal E}_m(T^*,V)<+\infty$. \label{inviscid} \end{theo}

Note that there are no restrictions on the viscosity parameters other than positiveness and sufficient regularity (${\cal C}^m$). It seems physically justified to ask $\mu(\rho)$ and $\lambda(\rho)$ to be increasing with the density, but, like the sign of the slip coefficient, the signs of $\mu'$ and $\lambda'$ do no intervene technically. 
Also, the results are shown for a barotropic pressure law, but we can extend them to any positive, increasing ${\cal C}^1$ pressure law.

Finally, we expect that our results are also valid in any domain of $\mathbb{R}^3$ with a ${\cal C}^{m'}$ boundary, for $m'$ large enough, locally characterised by equations of type $z=\psi(y)$. As shown in \cite{MR}, the differences are technical, and extra difficulties arise only because the normal vector to the boundary is no longer constant.
\newline

Our results are obtained by classical arguments once a uniform estimate is shown. The key bound is the following.
\begin{theo} \underline{Uniform energy bound.}

Let $m\geq 7$ and $M_0>0$. We assume that ${\cal N}_m(t,F)<+\infty$ for every $t>0$, and that the initial value of $U$ satisfies ${\cal N}_m(0,U) \leq M_0$. Then, there exist $\eps_0>0$, $T^*>0$ and a positive increasing function $Q: \rplus\rightarrow \rplus$, with $Q(z)\geq z$, such that, for $\eps\leq \eps_0$ and $0\leq t\leq T^*$,
\begin{equation} {\cal E}_m(t,U)+\trinorme{1,\infty,t}{\derp{z}u_3}^2 + \eps\int_0^t \norme{m}{\nabla u(s)}^2 + \norme{m-1}{\nabla^2 u_\tau(s)}^2~ds \leq Q(2M_0) \label{th13} \end{equation} \label{thEE} \end{theo}

Let us outline the proof of \thref{thEE}. It is proved by showing that the left-hand side of (\ref{th13}) is bounded by
$$ Q(M_0) + (t+\eps) Q({\cal E}_m(t,U)+{\cal N}_m(t,F)). $$
The energy function we consider, ${\cal E}_m$, contains $L^2$ in space norms with many derivatives, on the first line of (\ref{emtu}), and conormal-Lipschitz norms on the second line. The first type of term is dealt with by performing energy estimates, in which we will have to control the commutators between the conormal derivatives and the operators that appear in the equation. 
In particular, we will make use of the symmetrisable hyperbolic-parabolic structure of the compressible model to get the energy estimates on $U$. For the $L^\infty$ norms, we will widely use an anisotropic Sobolev embedding theorem (\thref{sobo}), which is the main contributor to the restriction on $m$. 
For the terms in the second line of (\ref{emtu}), a maximum principle will provide us with bounds on $\trinorme{1,\infty,t}{\derp{z}u\tang}$, while the Duhamel formula for the equation satisfied by $\derp{z}\rho$, which is obtained by combining $\eps\times(\ref{NSdiv})$ with the third component of (\ref{NS}), will give us the bounds for the $W^{1,\infty}_{co}$ norms on $\derp{z}\rho$. 
A bootstrap argument closes the proof.

Note that, in the context of \tworef{unifexist}{inviscid}{Theorems}, $M_0$ can be arbitrarily small, as ${\cal N}_m(0,U)=0$. We will be able to prove the energy estimate taking initial conditions into account; this allows one to extend our results to less regular force terms or different initial values, providing the compatibility conditions yield uniform bounds on the norms of $U$ at $t=0$.
\newline

\textbf{\underline{Organisation of the paper.}} In the next section, we prove \tworef{unifexist}{inviscid}{Theorems}, assuming \thref{thEE}. The remaining sections will all be dedicated to proving this uniform estimate. 
Starting with some important commutator estimates in section 3, we then proceed to prove the bound \textit{a priori}, looking at each component of $(U,\derp{z}U)$ separately, and getting the required conormal and $L^\infty$ bounds on each of them: $U$ in section 4, the normal derivative of $u\tang$ in section 5, $\derp{z}u_3$ in section 6, and the normal derivative of the density in section 7. 
We conclude the proof of \thref{thEE} in section 8.

\section{Proof of \tworef{unifexist}{inviscid}{Theorems}}

In this section, we assume the uniform estimate in \thref{thEE}. To obtain \thref{unifexist}, we can use a classical fixed-point iteration method to get existence of solutions for any fixed $\eps>0$ on a time interval $[0,T^\eps]$ depending on the viscosity (see \cite{WW}, section 4.1). A bootstrap argument with the uniform bounds then yields a uniform existence time.
\newline

We will further detail the proof of \thref{inviscid}. We first get local convergence by using a standard compactness argument.

\begin{propo} \underline{Conormal compact embedding theorem.}

Let $T>0$ and $(U_n)_{n\in\mathbb{N}}$ be a bounded sequence of $H^m_{co}([0,T]\times\Omega)$, such that the sequence $(\nabla U_n)_n$ is bounded in $H^{m-1}_{co}([0,T]\times\Omega)$. 
Then we can extract a sub-sequence $(U_{n_j})_{j\in\mathbb{N}}$ such that, for every $\alpha\in\mathbb{N}^4$ with $|\alpha|\leq m-1$, $Z^\alpha U_{n_j}$ converges in $L^2_{loc}([0,T]\times\Omega)$ - we will say that $(U_n)$ is \underline{locally compact} in $H^{m-1}_{co}([0,T]\times\Omega)$.
\label{rellich} \end{propo}

The bound on the gradient of $U_n$ is crucial here. In other contexts in which conormal Sobolev spaces have been used, such as in \cite{Go} and \cite{Sp}, one normal derivation costs two conormal derivations, thus their $H^m_{co}$ spaces locally and compactly embed in $H^{m-2}_{co}$. 
But in our energy estimates on the isentropic Navier-Stokes system, we find that one normal derivative can be controlled by the same number of conormal derivatives. The proof of \propref{rellich} is similar to that in \cite{Sp}.
\newline

For a given $\eps>0$, let $U^\eps=(\rho^\eps-1,u^\eps)$ be the solution to (\ref{NSdiv})-(\ref{NS})-(\ref{NPB})-(\ref{NBC})-(\ref{limcond}) with viscosity coefficient $\eps$, given by \thref{unifexist}. 
The energy estimate in \thref{thEE} tells us that the family $({\cal E}_m(T^*,U^{\eps}))_{0<\eps<\eps_{0}}$ is bounded, so we can immediately state that the sequence $(U^\eps)_{0<\eps<\eps_0}$ is locally compact in $H^{m-1}_{co}([0,T^*]\times \Omega)$ by \propref{rellich}. 
As $H^{m-1}_{co}([0,T^*]\times \Omega) \hookrightarrow {\cal C}([0,T^*],L^2(\Omega))$, we can consider a sequence $\eps_{n}\stackrel{n\rightarrow+\infty}{\rightarrow}0$ such that $U^{\eps_{n}}$ converges locally in $H^{m-1}_{co}([0,T^*]\times \Omega)$ and in ${\cal C}([0,T^*],L^2(\Omega))$ towards a function $V=(\rho-1,u)$, which is easily seen to be a weak solution to the compressible Euler system. 
Thanks to the uniform bounds, we see that $V$ has the same regularity as $U^{\eps}$, and in particular $V$ is Lipschitz-class, which yields uniqueness of solutions for the Euler equation in the space of functions satisfying ${\cal E}_m(T^*,V)<+\infty$. 
The whole sequence then converges towards $V$, strongly and locally in $H^{m-1}_{co}([0,T^*]\times\Omega)$, and $U^\eps(t)$ converges weakly in $L^2([0,T^*]\times\Omega)$ towards $V(t)$.
\newline

We now prove strong convergence in $L^2$. For $t\leq T^*$, we start with the classical energy inequality for the Navier-Stokes system, as in \cite{LM},
\begin{equation} \left| E(t,U^\eps) - E(0,U^\eps) - \int_0^t \int_\Omega \rho^\eps F(s)\cdot u^\eps (s) ~dx~ds \right| \leq C\eps\int_0^t \norme{H^1(\Omega)}{u^\eps(s)}^2~ds , \label{clasen} \end{equation}
$$ \mathrm{where}\hspace{10pt} E(t,U^\eps) = \int_\Omega \frac{1}{2}\rho^\eps(t)|u^\eps(t)|^2 + \frac{k}{\gamma(\gamma-1)}\int_\Omega (\rho^\eps)^\gamma(t) - 1 - \gamma(\rho^\eps(t)-1)~dx. $$
Note that $E(0,U^\eps) = 0$ for every $\eps$. As $\norme{H^1([0,T^*]\times\Omega)}{u^\eps}$ is bounded by \thref{thEE} and $\rho^\eps u^\eps$ converges weakly towards $\rho u$ in $L^2$ (as that is the case in the sense of distributions by local strong convergence), we get that, for $t\in [0,T^*]$,
$$ E(t,U^\eps) \stackrel{\eps\rightarrow 0}{\longrightarrow} \int_0^t \int_\Omega \rho(s,x) F(s,x)\cdot u(s,x)~dx~ds . $$
But $\int_{[0,t]\times\Omega} \rho F\cdot u$ is equal to $E(t,V)$, as (\ref{clasen}) with $\eps=0$ yields the energy equality for the compressible Euler equation. So, we have $E(t,U^\eps) \stackrel{\eps\rightarrow 0}{\longrightarrow} E(t,V)$ for every $t\leq T^*$.

We get $L^2$ convergence of the density by using the strict convexity of $z\mapsto z^\gamma$ for $z>0$ and $\gamma>1$. On one hand, we have
$$ \rho^\eps|u^\eps|^2-\rho|u|^2 - 2 \sqrt{\rho}u \cdot (\sqrt{\rho^\eps}u^\eps - \sqrt{\rho}u) = |\sqrt{\rho^\eps}u^\eps - \sqrt{\rho}u|^2 , $$
and on the other hand, there exists $c>0$ such that $\rho(t,x)>c>0$, thus the Taylor expansion of $z\mapsto z^\gamma$ yields
$$ (\rho^\eps)^\gamma - \rho^\gamma - \gamma \rho^{\gamma-1}(\rho^\eps-\rho) \geq c\gamma(\gamma-1) (\rho^\eps-\rho)^2 . $$
Adding the two together with the coefficients that appear in $E$, and taking the integral on $\Omega$, we get
\begin{equation} E(t,U^\eps)-E(t,V) - R(t,U^\eps,V) \geq c_0 \left(\norme{L^2(\Omega)}{(\sqrt{\rho^\eps}u^\eps - \sqrt{\rho}u)(t)}^2 + \norme{L^2(\Omega)}{(\rho^\eps-\rho)(t)}^2\right) \label{strconv} \end{equation}
where $R(t,U^\eps,V) = \int_\Omega \sqrt{\rho}u \cdot (\sqrt{\rho^\eps}u^\eps - \sqrt{\rho}u) + \frac{k}{\gamma-1}(\rho^\eps-\rho)(\rho^{\gamma-1}-1)~dx$, which converges to zero. 
Indeed, the uniform boundedness of the energy means that $\norme{L^2}{\sqrt{\rho^\eps(t)}u^\eps(t)}$ is bounded, so we can extract a weakly converging sub-sequence in $L^2(\Omega)$, $\sqrt{\rho^{\eps_n}(t)}u^{\eps_n}(t)$, whose limit is necessarily $\sqrt{\rho(t)} u(t)$, thus the whole sequence converges weakly in $L^2(\Omega)$. 
So the first term of $R(t,U^\eps,V)$ goes to zero. Likewise, we have $\rho^\eps(t) \rightharpoonup \rho(t)$ in $L^2(\Omega)$, and $\rho^{\gamma-1}-1$ is seen to be in $L^2(\Omega)$ by use of the order-one Taylor expansion of $z\mapsto z^{\gamma-1}$ at $z=1$, so the second term also converges to zero.

Moreover, $E(t,U^\eps)-E(t,V) \stackrel{\eps\rightarrow 0}{\longrightarrow} 0$, so (\ref{strconv}) gives us the global, strong $L^2$ convergence of $\rho^\eps$ towards $\rho$. Now, by remarking that
\begin{equation} \left|\norme{L^2(\Omega)}{(\sqrt{\rho^\eps}u^\eps)(t)}^2 - \norme{L^2(\Omega)}{(\sqrt{\rho}u^\eps)(t)}^2\right| \leq \norme{L^\infty}{u^\eps(t)}\norme{L^2}{u^\eps(t)}\norme{L^2}{\rho^\eps(t)-\rho(t)} , \label{normu} \end{equation}
we get that $\displaystyle \lim_{\eps\rightarrow 0} \norme{L^2(\Omega)}{(\sqrt{\rho} u^\eps)(t)} = \lim_{\eps\rightarrow 0} \norme{L^2(\Omega)}{(\sqrt{\rho^\eps}u^\eps)(t)} = \norme{L^2(\Omega)}{\sqrt{\rho}(t) u(t)}$, as $u^\eps$ is uniformly bounded in $L^2$ and $L^\infty$. 
So $u^\eps$ converges towards $u$ in $L^2([0,T^*]\times \Omega,dt~\rho dx)$, and $\rho dx$ is an equivalent measure to the Lebesgue measure, hence we conclude that $U^\eps$ converges towards $V$ in $L^2([0,T^*]\times\Omega)$. 
$L^\infty$ convergence is obtained by using the Sobolev embedding inequality and the uniform bounds. $\square$

\section{Preliminary properties of conormal derivatives}

We begin this section by reminding the reader of some important properties of our functional setting, and that will be used throughout the \textit{a priori} estimation process. Set $T>0$.

\begin{propo} \underline{Generalised Sobolev-Gagliardo-Nirenberg inequality (or tame estimate).} \cite{Go}

There exists a constant $C(T)$, which does not blow up as $T\rightarrow 0$, such that, for $f,g\in L^{\infty}([0,T]\times\Omega)\cap H^{m}_{co}([0,T]\times\Omega)$, and $\alpha_{1},~\alpha_{2}\in\mathbb{N}^{4}$ such that $|\alpha_{1}|+|\alpha_{2}|=m$,
$$ \int_0^T \norme{0}{(Z^{\alpha_{1}}f Z^{\alpha_{2}}g)(s)}^{2}~ds \leq C(T)\left[\trinorme{\infty,T}{f}^{2}\int_0^T \norme{m}{g(s)}^{2}~ds +\trinorme{\infty,T}{g}^{2}\int_0^T \norme{m}{f(s)}^{2}~ds \right].$$ \label{tame} \end{propo}

\begin{propo} \underline{Trace inequality.}

If $f\in L^2([0,T]\times\Omega)$ and $\nabla f\in L^2([0,T]\times\Omega)$, then
$$ \int_0^T \int_{\dOmega} |f(t,y,0)|^{2}~dt~dy \leq C \int_0^T \norme{0}{f(t)}\norme{0}{\derp{z}f(t)}~dt $$ \label{trace} \end{propo}

\begin{propo} \underline{Anisotropic Sobolev embedding theorem.}

If $f\in H^3_{co}([0,T]\times\Omega)$, $\nabla f\in H^2_{co}([0,T]\times\Omega)$, then $f\in L^\infty([0,T]\times\Omega)$ and
$$ \trinorme{\infty,T}{f}^{2} \leq C\left(\norme{2}{f(0)}^2+\norme{1}{\derp{z}f(0)}^2 + \int_0^T \norme{3}{f(t)}^{2}+\norme{2}{\derp{z}f(t)}^2 ~dt\right) $$ \label{sobo} \end{propo}

This last theorem is a direct application of the $H^m_{co}(\Omega)$ Sobolev embedding, used in \cite{MR} and \cite{MRfs}: for a given $t$, we have
$$ \norme{\infty}{f(t)}^2 \leq C(\norme{H^2_{co}(\Omega)}{f(t)}^2 + \norme{H^1_{co}(\Omega)}{\nabla f(t)}^2) \leq C(\norme{2}{f(t)}^2+\norme{1}{\derp{z}f(t)}^2). $$
We combine this with the following property: for $f\in H^{m+1}_{co}([0,T]\times\Omega)$,
\begin{equation} \trinorme{m,t}{f}^2 \leq \norme{m}{f(0)}^2 + C \int_0^t \norme{m+1}{f(s)}^2~ds . \label{l2time} \end{equation}
This is shown by writing $n(t)=n(0)+\int_0^t n'(s)~ds$ for $n(s)=\norme{m}{f(s)}^2$; given that
$$ n'(s) = 2\sum_{|\alpha|\leq m} \int_\Omega (\derp{t}Z^\alpha f(s,x))(Z^\alpha f(s,x))~dx, $$
we easily see that $n'(s)\leq 2\norme{m+1}{f(s)}^2$. 
Morally, the bound (\ref{l2time}) means that we can exchange an $L^\infty$-in-time norm for an $L^2$-in-time norm for the cost of one conormal derivative, similarly to (\ref{l2timeinf}).

We will now show the important commutator properties we will need.

\subsection{Commuting with $\derp{z}$}

The first technical key to the proof in the subsequent parts of the paper is how to estimate the commutators that will appear when applying $Z^{\alpha}$ to equations (\ref{NSdiv}), (\ref{NS}). 
A lot of the commutators are trivial, since all $Z_{i}$'s commute with $\derp{t}$, $\derp{y_{1}}$ and $\derp{y_{2}}$, but $Z_{3}$ does not commute with $\derp{z}$. 
Specifically, we have
\begin{equation} [Z_{3},\derp{z}] = \phi(z)\derp{zz}-\derp{z}(\phi(z)\derp{z}) = -\phi'(z)\derp{z} \label{commu1} \end{equation}
Likewise, we can observe the commutator of $Z_3$ with $\derp{zz}$, which will come from the $\div\Sigma$ term of (\ref{NS}). This time, we have
\begin{equation} [Z_{3},\derp{zz}] = -2\phi'\derp{zz} - \phi''\derp{z} . \label{commul} \end{equation}

When commutating with a higher order operator, $Z_3^m$ for $m>1$, we show the following:
\begin{propo} (a) For $m\geq 1$, there exist two families of bounded functions $(\varphi_{\beta,m})_{0\leq \beta <m}$ and $(\varphi^{\beta,m})_{0\leq\beta<m}$, such that
\begin{equation} [Z_{3}^{m},\derp{z}] = \sum_{\beta=0}^{m-1} \varphi_{\beta,m}(z)Z_{3}^{\beta}\derp{z} = \sum_{\beta=0}^{m-1} \varphi^{\beta,m}(z)\derp{z}Z_{3}^{\beta} \label{commum} \end{equation}

(b) For $m\geq 1$, there exist four families of bounded functions $(\psi_{1,\beta,m})$, $(\psi_{2,\beta,m})$, $(\psi^{1,\beta,m})$ and $(\psi^{2,\beta,m})$, for $0\leq \beta <m$, such that
\begin{eqnarray*} [Z_{3}^{m},\derp{zz}] & = & \sum_{\beta=0}^{m-1} \psi_{1,\beta,m}(z) Z_{3}^{\beta}\derp{z} + \psi_{2,\beta,m}(z)Z_{3}^{\beta}\derp{zz} \\
 & = & \sum_{\beta=0}^{m-1} \psi^{1,\beta,m}(z)\derp{z}Z_{3}^{\beta} + \psi^{2,\beta,m}(z)\derp{zz}Z_{3}^{\beta}. \end{eqnarray*}
\label{comprop} \end{propo}
In practice, we can therefore choose to place the normal derivatives as the first or last derivative to be applied in all the terms of the commutator. We can deduce from the proposition the basic estimate for the commutators with $\derp{z}$:
\begin{coro} For any $f\in H^{m}_{co}([0,T]\times \Omega)$ such that $\derp{z}f\in H^{m-1}_{co}([0,T]\times \Omega)$, and for any $|\alpha|\leq m$,
$$ \int_0^T \norme{0}{[Z^\alpha,\derp{z}]f(t)}^2~dt \leq C\int_0^T \norme{m-1}{\derp{z}f(t)}^2~dt $$ \label{comcoro} \end{coro}
We will deal with the commutators with $\derp{zz}$ directly in context: they often appear with a factor $\eps$, and equation (\ref{NS}) will allow us to substitute the difficult terms.
\newline

\preu{Proof of \propref{comprop}}: equation (\ref{commu1}) shows the case $m=1$, and we continue by induction. Let us just explain the case $m=2$ to show the mechanism; the rest is left to the reader. We have
\begin{eqnarray*} [Z_{3}^{2},\derp{z}] & = & Z_{3}[Z_{3},\derp{z}]+[Z_{3},\derp{z}]Z_{3} \\
 & = & -\phi\phi''\derp{z}-\phi'(Z_{3}\derp{z}+\derp{z}Z_{3}). \end{eqnarray*}
Using (\ref{commu1}), we can write the second part of the last line as either $\phi'(z)(2Z_{3}\derp{z}+\phi'\derp{z})$ or $\phi'(z)(2\derp{z}Z_{3}-\phi'\derp{z})$, which proves the proposition for $m=2$.

The proof of (b) is also an elementary induction. $\square$

\subsection{Commuting with a function}

\begin{nota} For $\alpha,~\beta\in\mathbb{N}^{4}$, we write $\beta \leq \alpha$ if, for every $i$, $\beta_{i}\leq\alpha_{i}$. \end{nota}

\begin{propo} Let $\alpha\in\mathbb{N}^{4}$, $|\alpha|=m>0$, be fixed, and $f\in H^m_{co}([0,T]\times\Omega)\cap L^\infty([0,T]\times\Omega)$ and $g\in L^\infty$ such that $\derp{t}g,~\nabla g\in H^{m-1}_{co}([0,T]\times\Omega)\cap L^\infty([0,T]\times\Omega)$. Then we have the following inequality:
\begin{eqnarray} \int_{0}^{T}\norme{0}{[Z^{\alpha},g]f}^{2}~dt & \leq & C\sum_{j=0}^3\int_{0}^{T}\trinorme{\infty,T}{Z_j g}^{2}\norme{m-1}{f}^{2}+\trinorme{\infty,T}{f}^{2}\norme{m-1}{Z_{j}g}^{2}~dt, \label{funcom} \\ & \leq & C\int_{0}^{T}\trinorme{1,\infty,T}{g}^{2}\norme{m-1}{f(s)}^{2}+\trinorme{\infty,T}{f}^{2}\norme{m}{g(s)}^{2}~dt, \nonumber \end{eqnarray}
if, moreover, $g\in L^2([0,T]\times\Omega)$. \label{funcprop} \end{propo}

This proposition is easily proved, using the Leibniz formula and \propref{tame}. We will require formulation (\ref{funcom}) whenever $g\notin L^2$ (for example when $g=\rho$).
\newline

We prove one more commutator estimate, which we will need when estimating the normal derivatives of $u_\tau$ in the conormal spaces (section 5), as directly using the above would lead to $\norme{m-2}{\derp{zz} u}$ appearing, which we cannot bound uniformly in $\eps$ by using the equation.
\begin{propo} Let $f,~g,~\alpha$ be as in \propref{funcprop}, with $g$ scalar such that $g|_{z=0}=0$. Then there exists $C$ which does not depend on $T$ such that
$$ \int_{0}^{T}\norme{0}{[Z^{\alpha},g\derp{z}]f}^{2}~dt \leq C\int_{0}^{T}(\trinorme{\Lip,T}{g}^{2}+\trinorme{1,\infty,T}{\nabla g}^{2})\norme{m}{f(s)}^{2}+\trinorme{1,\infty,T}{f}^{2}\norme{m-1}{\nabla g(s)}^{2}~dt $$ \label{comlem} \end{propo}

\preu{Proof:} we decompose the commutator as $[Z^{\alpha},g\cdot\nabla]f = [Z^{\alpha},g]\derp{z} f + g [Z^{\alpha},\derp{z}]f$, and start by taking a closer look at the second term. 
By \propref{comprop}, it is equal to a sum of terms of the form $\varphi^{\beta}(z)g~\derp{z}Z^{\beta}f$, with $\beta\leq\alpha$, $\beta\neq\alpha$. As $g|_{z=0}=0$, we have
\begin{equation} |g(t,x)|\leq \phi(z)\trinorme{\infty,T}{\derp{z}g}, \label{u3phi} \end{equation}
so, as $Z_{3}=\phi(z)\derp{z}$,
\begin{equation} \norme{0}{g[Z^{\alpha},\derp{z}]f} \leq C\trinorme{\infty,T}{\nabla g}\norme{m-1}{\phi(z)\derp{z}f} \leq C\trinorme{\Lip,T}{g}\norme{m}{f}. \label{hmcom3} \end{equation}

Now we look at the first term of the decomposition. We prove that
\begin{equation} \int_{0}^{T}\norme{0}{[Z^{\alpha},g]\derp{z}f}^{2} \leq C\int_{0}^{T}(\trinorme{\Lip,T}{g}^{2}+\trinorme{1,\infty,T}{\nabla g}^{2})\norme{m}{f}^{2} + \trinorme{1,\infty,T}{f}^{2}\norme{m-1}{\nabla g}^{2}~dt, \label{comfinal} \end{equation}
which would end the proof of the proposition.

We can write $[Z^{\alpha},g]\derp{z}f$ as a sum, on $\beta$ and $\delta$, of terms of the form $\varphi^{\delta}(z) Z^{\beta}g~\derp{z}Z^{\delta} f$, where $\beta\leq\alpha$ and $|\beta|>0$, and $\delta\leq(\alpha-\beta)$. The idea is to insert $\frac{1}{\phi(z)}\times\phi(z)$, thus we have to estimate
$$ \norme{0}{\frac{1}{\phi(z)}(Z^{\beta}g)(Z_{3}Z^{\delta}f)}. $$
This will obviously be done by using the tame estimate, but we also have to deal with the first factor. We have $\phi^{-1}Z^{\beta}g = Z^{\beta}(\phi^{-1}g) - [Z^{\beta},\phi^{-1}]g$, and we need to write this last term more explicitly.
\begin{lemma} For every $b\in\mathbb{N}$, $\phi Z_{3}^{b}(\phi^{-1})$ is smooth and bounded, then there exists a family of smooth bounded functions $(\sigma_{p})_{0\leq p\leq b}$ such that for every $f\in H^{m}_{co}([0,T]\times\Omega)$,
$$ [Z^{m}_{3},\phi^{-1}]f = \sum_{p=0}^{m} \sigma_{p}(z) Z^{p}_{3}(\phi^{-1}f) $$ \label{phicom} \end{lemma}

\preu{Proof:} as $\phi$ only depends on $z$, we are only interested in $[Z_3^m,\phi^{-1}]$ for $m\geq 1$. The Leibniz formula yields that the commutator is a linear combination of terms written as $(Z_3^b(\phi^{-1}))(Z^{m-b}f)$, for $0< b\leq m$. 
In the case of $\phi(z)=\frac{z}{1+z}$, we notice that $Z_3^b (\phi^{-1})$ has the same following properties as $\phi^{-1}$: it is smooth on $]0,+\infty[$, bounded at infinity and has the same blow-up rate at $z=0$ as $\phi^{-1}$ (blows up like $z^{-1}$). So, for each $b$, $\phi Z^b_3(\phi^{-1})$ is a bounded function on $[0,+\infty[$, and if we write
$$ (Z^b_{3}(\phi^{-1}))(Z^{m-b}f) = \phi Z^b_{3}(\phi^{-1})\left[Z^{m-b}\left(\frac{1}{\phi}f\right) - [Z^{m-b},\phi^{-1}]f\right], $$
we have $m-b<m$, and we can reiterate the process. As $m$ is fixed, we conclude the proof with a finite number of iterations. $\square$
\newline

Applying the lemma, we have that, in any norm,
\begin{equation} \norme{}{\phi^{-1}Z^{\beta}g} \leq C\sum_{\beta'\leq\beta} \norme{}{Z^{\beta'}(\phi^{-1}g)} \label{phicom2} \end{equation}
With that, the tame estimate gives us
$$ \int_{0}^{T}\norme{0}{\phi^{-1}(Z^{\beta}g)(Z_{3}Z^{\delta}f)}^{2}~dt \leq C \int_{0}^{T} \trinorme{\infty,T}{Z_{3}f}^{2} \sum_{j=0}^{3} \norme{m-2}{Z_j\left(\frac{g(s)}{\phi}\right)}^{2} $$
$$ \hspace{70pt} + \left(\trinorme{\infty,T}{\phi^{-1}g}^2 +\sum_{j=0}^{3} \trinorme{\infty,T}{Z_{j}(\phi^{-1}g)}^2\right)\norme{m-1}{Z_3 f(s)}^{2}~dt $$
and the $\trinorme{\infty,T}{\phi^{-1}g_3}$ comes from the terms in (\ref{phicom2}) with $\beta'=0$.

It remains to deal with the terms involving $\phi^{-1}g$, and the key fact here is that $Z^{\beta}g|_{z=0}=0$. We easily have by (\ref{u3phi}),
\begin{equation} \norme{\infty}{\phi^{-1}g} \leq \norme{\Lip}{g} \mathand \norme{\infty}{Z_{j}(\phi^{-1}g)}\leq C\norme{1,\infty}{\nabla g}, \label{hardyinf} \end{equation}
as $Z_{3}(\phi^{-1}g) = \phi^{-1}Z_{3}g+\phi'\phi^{-1}g$, with $\phi'$ bounded, and likewise we can use the Hardy inequality to get
$$ \norme{m-2}{Z_{j}(\phi^{-1}g)} \leq C\norme{m-1}{\derp{z}g}, $$
which ends the proof of (\ref{comfinal}). $\square$

\section[Proof of \thref{thEE}, part I: \textit{a priori} estimates on $U$]{Proof of \thref{thEE}, part I \\ \textit{A priori} estimates on $U$
 \sectionmark{Proof of \thref{thEE} part I: a priori estimates on $U$}}
\sectionmark{Proof of \thref{thEE} part I: a priori estimates on $U$}

\begin{assu} Throughout the \textit{a priori} estimation process, we will assume that there is no vacuum on the time of study: there exists $0<c_{0}<1$ such that $c_{0}\leq\rho(t,x)$ for $x\in\Omega$. 
Also, as $u_3|_{z=0}=0$, for any $\delta>0$, we assume that there exists $z_\delta>0$, independent of $\eps$, such that $|u_3(t,x)|\leq \delta$ for $x\in \sdel{F} := \mathbb{R}^2\times [0,\sdel{z}]$.

After getting the estimates, we will show in section 8 that the final bounds actually imply these two properties on $[0,T^*]$, with $T^*$ such that ${\cal E}_m(T^*,U)\leq M$, and prove by a bootstrap argument that $T^*$ does not depend on $\eps$.
\label{assu1} \end{assu}

\begin{nota} As of now, $0<c<1$ will designate a small constant, $C>1$ a large constant, and $Q(z)$ a positive increasing function of $\rplus$, with polynomial growth and $Q(z)\geq z$. 
All three can change from one line to the next, and can depend on any of the system's parameters (constants $a$, $k$ and $\gamma$, or bounds of the viscosity functions $\mu(\rho)$ and $\lambda(\rho)$ - we will often omit the dependence on $\rho$), on the order of derivation $m$ or on $\eps_0$. \end{nota}

\subsection{Conormal energy estimates}

We start with energy estimates on $U=(\rho-1,u)$ and its conormal derivatives: this will estimate the $\trinorme{m}{U}$ part of ${\cal E}_{m}(t,U)$, and it will uncover other terms of ${\cal E}_{m}(t,U)$ that we will need to estimate later.

We will be able to estimate $\rho-1$ and $u$ simultaneously by taking full advantage of the symmetrisable hyperbolic structure of the order-one part of the isentropic Navier-Stokes system $(\ref{NSdiv}),~ (\ref{NS})$. We rewrite it as
\begin{equation} A_{0}(\rho)\derp{t}U + \sum_{j=1}^{3} A_{j}(U)\derp{x_{j}}U - \eps(0,\mu(\rho) \Delta u + \lambda(\rho) \nabla\div u)^t = (0,\rho F-\eps\sigma(\nabla U))^t , \label{matrixform} \end{equation}
where $\lambda=\mu+\lambda_0>0$, $\sigma$ has the following expression
$$ \sigma(\nabla \rho,\nabla u) = 2Su\cdot\nabla (\mu(\rho)) + \div u \nabla (\lambda_0(\rho)) = 2\mu'(\rho) Su\cdot \nabla \rho + \lambda'(\rho)\div u \nabla \rho. $$
The matrices $A_{j}(U)$ are: $A_{0}(\rho)=\diag(1,\rho,\rho,\rho)$,
$$ A_{1}(U)=\left(\begin{array}{cccc} u_{1} & \rho & 0 & 0 \\ k\rho^{\gamma-1} & \rho u_{1} & 0 & 0 \\ 0 & 0 & \rho u_{1} & 0 \\ 0 & 0 & 0 & \rho u_{1} \end{array} \right) ~,~ A_{2}(U)=\left(\begin{array}{cccc} u_{2} & 0 & \rho & 0 \\ 0 & \rho u_{2} & 0 & 0 \\ k\rho^{\gamma-1} & 0 & \rho u_{2} & 0 \\ 0 & 0 & 0 & \rho u_{2} \end{array} \right) $$
$$ \mathrm{and} ~ A_{3}(U)=\left(\begin{array}{cccc} u_{3} & 0 & 0 & \rho \\ 0 & \rho u_{3} & 0 & 0 \\ 0 & 0 & \rho u_{3} & 0 \\ k\rho^{\gamma-1} & 0 & 0 & \rho u_{3} \end{array} \right). $$
This system is symmetrisable: multiplying these matrices on the left by the positive diagonal matrix $D(\rho)=\diag(k\rho^{\gamma-2},1,1,1)$, we get:
\begin{equation} DA_{0}\derp{t}U + \sum_{j=1}^{3} DA_{j}\derp{x_{j}}U -\eps(0,\mu\Delta u +\lambda\nabla\div u)^t = (0,\rho F-\eps\sigma(\nabla U))^t \label{sym} \end{equation}

$DA_{0}(\rho)$ is symmetric, so
$$ \frac{1}{2} \frac{d}{dt}\left(\int_{\Omega} DA_{0}U\cdot U\right) = \int_{\Omega} DA_{0}\derp{t}U\cdot U + \frac{1}{2}\int_{\Omega} \derp{t}(DA_{0})U\cdot U. $$
Integrating this in time between $0$ and $t$, and since $\rho$ is uniformly bounded from below by $c_{0}$ and from above by $c_{1}$, there exist $0<c<C$ such that
\begin{equation} c \norme{0}{U(t)}^{2} \leq C\norme{0}{U(0)}^{2}+C\trinorme{\infty,t}{\derp{t}DA_0}\int_{0}^{t}\norme{0}{U(s)}^{2} + \int_{0}^{t}\int_{\Omega} DA_{0}\derp{t}U(s)\cdot U(s)~ds \label{l2ee1} \end{equation}

We replace $DA_{0}\derp{t}U$ by its expression in (\ref{sym}), and we use integration by parts on the integrals with order-one derivatives of $U$:
\begin{equation} \int_{\Omega} DA_{j}\derp{x_{j}}U \cdot U = -\frac{1}{2} \int_{\Omega} \derp{x_{j}}(DA_{j}) U\cdot U . \label{l2ee2} \end{equation}
Indeed, $DA_{j}(U)$ is symmetric, so $DA_{j}U\cdot \derp{x_{j}}U = U\cdot DA_{j}\derp{x_{j}}U$, and we notice that, for each $j$, $DA_{j}U\cdot U = (3\rho P'(\rho)+\rho|u|^2)u_j$, which means that there is no boundary term (at $z=0$ for $j=3$) when we integrate by parts, which gives us
\begin{eqnarray} c\norme{0}{U(t)}^2 & \leq & C\norme{0}{U(0)}^2 + C\left(\sum_{j=0}^3 \trinorme{\Lip,t}{DA_j}\right)\int_0^t \norme{0}{U(s)}^2 ~ds \nonumber \\
 & & \hspace{10pt} + \eps\int_0^t \int_\Omega \mu\Delta u\cdot u + \lambda \nabla \div u \cdot u + \sigma\cdot u~ds + \int_0^t \int_\Omega \rho F\cdot u~ds \nonumber \\
 & \leq & C\norme{0}{U(0)}^2 + C(1+\trinorme{\Lip,t}{U}^2)\int_0^t \norme{0}{U(s)}^2 ~ ds \nonumber \\
 & & \hspace{10pt} + \eps \int_0^t \int_\Omega \mu\Delta u \cdot u + \lambda \nabla\div u\cdot u + \sigma\cdot u ~ds + \int_0^t \int_\Omega \rho F\cdot u~ds \label{l2ee2b} \end{eqnarray}

We use integration by parts on the order-two derivatives and the Navier boundary condition (\ref{NBC}) to deal with the boundary term:
\begin{eqnarray} \int_\Omega \mu\Delta u\cdot u + \lambda \nabla\div u\cdot u & = & -\int_\Omega (\mu|\nabla u|^2 + \lambda|\div u|^2) -\int_{z=0} 2a |u\tang|^2 \hspace{50pt} \nonumber \\
 & & \hspace{20pt} - \int_\Omega \left(\nabla (\mu(\rho)) \cdot \nabla u\cdot u + \div u~ \nabla (\lambda(\rho))\cdot u\right) , \label{l2ee3} \end{eqnarray}
with the notation $\int_{z=0}f = \int_{\rtwo} f(y,0)~dy$. The first term of (\ref{l2ee3}) is moved the left-hand side of (\ref{l2ee2b}), as is the second if $a>0$. 
When $a<0$, we can use the trace theorem, \propref{trace}, and absorb the norm of $\nabla u$ in the left-hand side by using Young's inequality, $\nu\zeta \leq \frac{\eta}{2} \nu^{2} + \frac{1}{2\eta} \zeta^{2}$ for any $(\nu,\zeta)\in\mathbb{R}^2$ and $\eta>0$, with an adequate parameter $\eta$ (same as in \cite{Pm11}). 
Young's inequality is also used on the term containing the derivatives of $\mu$ and $\lambda$, which turns (\ref{l2ee3}) into
\begin{equation} \int_{\Omega} (\mu\Delta u\cdot u+\lambda\nabla \div u\cdot u) \leq \frac{c(1+\trinorme{\Lip,t}{\rho})}{\eta}\norme{0}{u}^2 - (c-(|a|+1)\eta) \norme{0}{\nabla u}^2 , \label{l2ee3b} \end{equation}
so we choose $\eta$ so that $c-(|a|+1)\eta = c/2$, which allows us to move the $\norme{0}{\nabla u}^2$ term to the left-hand side and absorb it with the first term of (\ref{l2ee3}) that we moved there earlier. The remainder $\int_\Omega \sigma\cdot u$ is bounded the same way.

Combining (\ref{l2ee2b}) and (\ref{l2ee3b}), and the assumption that $\rho$ is uniformly bounded, there exist $0<c<C$ such that
$$ c\left[\norme{0}{U(t)}^{2} + \eps\int_{0}^{t}\norme{0}{\nabla u}^{2} \right] \leq C\left(\norme{0}{U(0)}^{2}+(1+\trinorme{\Lip,t}{U}^2) \int_{0}^{t}\norme{0}{U(s)}^{2} + \norme{0}{F(s)}^2~ds\right), $$
as $\norme{\infty}{\nabla(DA_{j})}\leq C\norme{\infty}{U}\norme{\infty}{\nabla U}$.
\newline

We have shall now show the following, higher order estimate.

\begin{esti} For every $m\geq 0$,
$$ c\left[\norme{m}{U(t)}^{2} + \eps\int_{0}^{t}\norme{m}{\nabla u(s)}^{2} ~ds\right] \hspace{250pt} $$
$$ \hspace{10pt} \leq C\left[\norme{m}{U(0)}^{2}+(1+\trinorme{\Lip,t}{U}^2+\trinorme{\infty}{F}^2)\int_{0}^{t} \norme{m}{U}^{2} + \norme{m-1}{\nabla U}^{2} + \norme{m}{F}^2 ~ds\right] $$
(the negative index that appears when $m=0$ is ignored). \label{hmee} \end{esti}

\preu{Proof:} higher-order estimates work exactly the same way as above, only we will have to estimate the commutators between $Z^{\alpha}$, with $|\alpha|\leq m$, and the operators in the equation. 
We apply $Z^{\alpha}$ to equation (\ref{matrixform}), and isolate the highest-order terms as follows:
\begin{equation} A_0 \derp{t}Z^{\alpha}U + \sum_{j=1}^{3} A_{j}\derp{x_{j}}(Z^{\alpha}U) - \eps(0,\div Z^\alpha(\mu\nabla u)+\lambda\nabla\div (Z^{\alpha}u))^t = (0,Z^\alpha (\rho F+\eps\tilde{\sigma}))^t-{\cal C}^{\alpha}. \label{matrixformalpha} \end{equation}
${\cal C}^{\alpha}$ contains the commutators:
$$ {\cal C}^{\alpha} = [Z^\alpha,A_0\derp{t}]U + \sum_{j=1}^{3} [Z^{\alpha},A_{j}\derp{x_{j}}]U -\eps(0, [Z^{\alpha},\div] (\mu \nabla u)+[Z^{\alpha},\lambda(\rho)\nabla \div]u)^t. $$
Notice the peculiar formulation for the laplacian term. 
This allows us to avoid difficult commutator terms on the boundary when $\mu$ is not constant; indeed, the Navier boundary condition for $Z^\alpha u$ is $Z^\alpha(\mu\derp{z}u\tang) = 2a Z^\alpha u\tang$, thus, when we multiply equation (\ref{matrixformalpha}) by $Z^\alpha u$ and integrate by parts, we can use this boundary condition immediately. 
The remainder of the order-two part of the equation $\tilde{\sigma}$ is the equivalent of $\sigma$ for this formulation.

We then multiply by the matrix $D$, which is uniformly bounded in $L^\infty$, and can repeat the above to obtain:
$$ c\left[\norme{0}{Z^{\alpha}U(t)}^{2} + \eps\int_{0}^{t}\left(\norme{0}{Z^\alpha \nabla u}^{2} + \int_{z=0}|Z^{\alpha}u|^{2}\right)\right] \leq C\norme{0}{Z^{\alpha}U(0)}^{2} \hspace{100pt} $$
$$ \hspace{80pt} +C\int_0^t (1+\trinorme{\Lip,t}{U}^2+\trinorme{\infty,t}{F})\norme{m}{U}^{2} + \eps \norme{m-1}{\nabla U}^2 + \norme{m}{F}^2 ~ds $$
\begin{equation} \hspace{120pt} + C\int_0^t\int_{\Omega} |(\eps Z^\alpha \sigma+{\cal C}^{\alpha})\cdot Z^{\alpha}U| + |{\cal I}^\alpha_\nabla| ~ds. \label{hmee1} \end{equation}
The term ${\cal I}^\alpha_\nabla$ contains extra commutators arising from the following integration by parts:
\begin{eqnarray*} -\int_\Omega \div Z^\alpha(\mu\nabla u)\cdot Z^\alpha u & = & \int_{\dOmega} Z^\alpha(\mu\derp{z}u)\cdot Z^\alpha u + \int_\Omega Z^\alpha(\mu \nabla u):\nabla Z^\alpha u \\
 & = & \int_{\dOmega} 2a |Z^\alpha u\tang|^2 + \int_\Omega \mu|Z^\alpha \nabla u|^2 \\
 & & \hspace{40pt} +\int_\Omega ([Z^\alpha,\mu]\nabla u) : Z^\alpha \nabla u - Z^\alpha(\mu\nabla u):([Z^\alpha,\nabla] u), \end{eqnarray*}
so ${\cal I}^\alpha_\nabla$ is the final integral multiplied by $\eps$. Here, we have used the contracted matrix product: $A:B = \sum_{i,j} a_{i,j} b_{i,j}$.

To deal with the commutator terms, we will use the tools shown in section 3. Let us first look at ${\cal I}^\alpha_\nabla$. In the first term, we have
$$ \eps \left| \int_\Omega ([Z^\alpha,\mu]\nabla u) : Z^\alpha \nabla u \right| \leq \norme{0}{[Z^\alpha,\mu]\nabla u} \times \eps\norme{0}{Z^\alpha \nabla u}. $$
By estimate (\ref{funcom}) in \propref{funcprop}, the integral in time of the square of the first norm is easily bounded, whereas the second can be moved to the left-hand side (absorbed) by using Young's inequality with an adequate parameter $\eta$ to have
\begin{equation} \left|\int_0^t \eps \norme{m}{\nabla u} \norme{0}{[Z^\alpha,\mu]\nabla u}\right | \leq C\eps (1+\trinorme{\Lip,t}{U}^2)\int_0^t (\norme{m}{\rho-1}^2+\norme{m-1}{\nabla u}^{2}) + \frac{c}{4}\eps\int_0^t \norme{0}{Z^\alpha \nabla u}^{2} . \label{hmeeb} \end{equation}
The second term of ${\cal I}^\alpha_\nabla$ is dealt with in exactly the same fashion, with \cororef{comcoro} controlling the side of the product containing the commutator, and Young's inequality allowing to absorb the other.
\newline

Now to the commutators in ${\cal C}^\alpha$ relative to the order-one part of the equation. We notice that, for $j\in\{0,1,2,3\}$,
$$ [Z^{\alpha},A_{j}\derp{x_j}]U = [Z^{\alpha},A_{j}]\derp{x_j} U+A_{j}[Z^{\alpha},\derp{x_j}]U, $$
where $x_0$ is understood to be $t$. The first term is estimated using inequality (\ref{funcom}) of \propref{funcprop} (we cannot write $\norme{m}{A_{j}}$ because $\rho \notin L^{2}$), and the second, which is either 0, if $j\neq 3$, or $A_{3}[Z^{\alpha},\derp{z}]U$, is estimated using \cororef{comcoro}, so estimating this commutator yields
\begin{equation} \int_{\Omega} |A_{3}[Z^{\alpha},\derp{z}]U\cdot Z^{\alpha}U | \leq C(1+\trinorme{\Lip,t}{U}^{2})\norme{m-1}{\nabla U}^{2}+\norme{m}{U}^{2}. \label{hmee2} \end{equation}

The commutators on the order-two part of the equation can be split into terms of two types as follows: commutators on the viscosity parameters (functions of $\rho$) and commutators on the differential operators, for instance
\begin{equation} [Z^\alpha,\lambda\nabla\div] = [Z^\alpha,\lambda]\nabla \div + \lambda[Z^\alpha,\nabla\div] . \label{lapcomsplit} \end{equation}

We begin with the first term. Notice that we can write this commutator as a sum of the following type of integral,
$$ I_{i,j} = \int_0^t \int_\Omega \eps [Z^\alpha,\lambda]\derp{x_i,x_j}u \cdot Z^\alpha u, $$
with $i,~j \in \{1,2,3\}$. Using \tworef{funcprop}{comprop}{Propositions}, we write
\begin{equation} I_{i,j} = \eps \sum_{|\beta|+|\delta|\leq m, |\delta|< m} \int_0^t \int_\Omega \psi_{\beta,\delta}(z) Z^\beta(\lambda(\rho))  ~ \derp{x_i}Z^\delta\derp{x_j}u\cdot Z^\alpha u ~ ds , \label{801} \end{equation}
where the $\psi_{\beta,\delta}$ are ${\cal C}^\infty$ bounded functions of $z$. If $|\beta|< m$, we can integrate by parts on the $x_i$ variable:
\begin{equation} I_{i,j} \leq \sum_{|\beta|+|\delta|\leq m} c_{\beta,\delta} \eps \int_0^t \int_\Omega |\derp{x_i}Z^\beta(\lambda(\rho))~Z^\delta\derp{x_j}u\cdot Z^\alpha u + Z^\beta(\lambda(\rho))~Z^\delta\derp{x_j}u\cdot \derp{x_i}Z^\alpha u|~ds . \label{802} \end{equation}
We use Young's inequality, then we use the tame estimate on the left of the scalar products to get
$$ I_{i,j} \leq \frac{c\eps}{2}\int_0^t \norme{m}{\nabla u}^2~ds + \eps Q(\trinorme{\Lip,t}{\rho}^2) \int_0^t \norme{m-1}{\nabla u}^2 + \norme{m}{u}^2~ds \hspace{50pt} $$
$$ \hspace{100pt} + \eps\trinorme{\infty,t}{\nabla u}^2 Q\left(\int_0^t \norme{m-1}{\nabla \rho}^2~ds\right) . $$
This allows us to absorb the first term. Note that we have used the tame estimate on the pair $(Z \rho,\derp{x_j}u)$ rather than $(\rho,\derp{x_j}u)$ on the second term of (\ref{802}) to get $H^{m-1}_{co}$ norms of $\nabla u$ with a factor $\eps$, which can be controlled by induction on $m$.

There is a slight difference when $\beta=\alpha$ in (\ref{801}). In this case, $\delta=0$, but after integrating by parts, we would need to control $\norme{m}{\nabla \rho}$. 
So we cannot integrate by parts and must estimate $J_{i,j}:=\eps\int_0^t \int_\Omega Z^\alpha(\lambda(\rho)) ~\derp{x_i,x_j}u\cdot u$ directly. 
This is not a problem if $(i,j)\neq (3,3)$: we bound $J_{i,j}$ by
$$ J_{i,j}\leq \eps \int_0^t \trinorme{\infty}{u}^2 \norme{1}{\nabla u}^2 Q\left(\norme{m}{\rho-1}^2\right)~ds . $$
At this stage, if $m=1$, we use Young's inequality to absorb $\eps\int_0^t \norme{1}{\nabla u}^2$, and if $m\geq 2$, we can use \propref{hmee} with $m=1$. When $(i,j)=(3,3)$, we instead replace $\eps\derp{zz}u$ by using the equation. We have
$$ \eps(\mu \derp{zz}u_\tau, (\mu+\lambda)\derp{zz}u_3) = \rho\derp{t}u + (\rho u\cdot\nabla u) + \nabla P - \rho F -\eps v_2, $$
where order-two derivatives appear in $v_2$. The terms that appear are therefore either controlled in $L^\infty$ by $\norme{\Lip}{U}$, or, in the case of $v_2$, dealt with using Young's inequality as above. 
The term $\int_\Omega \eps Z^\alpha\sigma\cdot Z^\alpha u$ in (\ref{hmee1}) is dealt with the same techniques as these commutators (the terms it contains ressemble those in (\ref{802})).
\newline

For the remaining commutator in (\ref{lapcomsplit}) (the estimation of the last term $[Z^\alpha,\div] (\mu\nabla u)$ follows the same lines), we assume that $\alpha_3>0$, as these commutators are zero otherwise. 
The terms of the commutator $[Z^{\alpha},\nabla\div]u$ of the form $[Z^{\alpha},\nabla\tang\derp{z}]u_{3}$ (the others are either trivial or contain $\derp{zz}$), can then be estimated with \cororef{comcoro}:
$$ \int_{\Omega} \eps\lambda|[Z^{\alpha},\derp{z}]\nabla\tang u_{3}\cdot Z^{\alpha}u\tang| \leq C\eps\norme{m}{\derp{z}u}\norme{m}{u} \leq C\eps\norme{m}{u}^{2}+\frac{c}{4}\eps\norme{m}{\derp{z}u}^{2}, $$
again using Young's inequality to allow us to absorb $\eps\int_0^t \norme{m}{\derp{z}u}^{2}$. On to the remaining term $\eps[Z^{\alpha},\derp{zz}]u$. By \propref{comprop}, there exist two families of functions, $(\varphi_{\beta})_{\beta\leq\alpha}$ and $(\psi^\beta)_{\beta\leq\alpha}$, such that
$$ [Z^{\alpha},\derp{zz}]u = \sum_{\beta\leq\alpha,~\beta\neq\alpha} \varphi_\beta(z) Z^{\beta}\derp{z}u+\psi^\beta(z) \derp{zz}Z^{\beta}u $$
The first part of the right-hand side is simply bounded by $\norme{m-1}{\derp{z}u}$, so it only remains to control $J:=\eps\int_{\Omega} |\derp{zz}Z^{\beta}u\cdot Z^{\alpha}u|$. As $\alpha_{3}>0$, we have $Z^{\alpha}u=0$ on the boundary, so integrating by parts and again using Young's inequality provides us with
$$ J = \eps\int_{\Omega} |\derp{z}Z^{\beta}u\cdot\derp{z}Z^{\alpha}u| \leq C\eps\norme{m-1}{\derp{z}u}\norme{m}{\derp{z}u} \leq C\eps\norme{m-1}{\derp{z}u}^{2}+\frac{c}{4}\eps\norme{m}{\derp{z}u}^{2}. $$

Combining this last inequality with (\ref{hmee1}), (\ref{hmeeb}) and (\ref{hmee2}), then summing for all $|\alpha|\leq m$, we get the result. $\square$

\subsection{$L^\infty$ estimates}

With conormal energy estimates, we obtain inequalities for any order of derivation $m$, and these inequalities contain $L^\infty$ or Lipschitz norms of $U$, and, as we shall see soon, $W^{1,\infty}_{co}$ norms of $\nabla U$. 
These need to be controlled, and the goal of this control is to get the conormal energy estimates to be closed for $m$ large enough, by using the Sobolev embedding inequality in \propref{sobo}. 
Putting the $L^\infty$ norm of the normal derivative $\derp{z}U$ to one side for now (this will be dealt with in the next sections), we can apply \propref{sobo} to $\trinorme{1,\infty,t}{U}$, and (\ref{l2time}) again on the terms involving the derivatives of $\rho$ and $u_3$.

\begin{esti} We have the following bound for $\trinorme{1,\infty,t}{U}$:
\begin{eqnarray*} \trinorme{1,\infty,t}{U}^2 & \leq & C(\norme{4}{U(0)}^2+\norme{4}{\derp{z}U(0)}^2) \\
 & & \hspace{15pt} + Ct\left(\trinorme{5,t}{U}^2+\trinorme{4,t}{\derp{z}u\tang}^2+\int_0^t \norme{5}{\derp{z}u_3}^2+\norme{5}{\derp{z}\rho}^2~ds\right) \end{eqnarray*} \label{linfu} \end{esti}

We can inject this into \propref{hmee} and partially close the energy estimate for $m-1\geq 5$. We widely use the Young inequality to separate product terms, and we get the following.

\begin{esti} For every $m\geq 6$, there exists a constant $C>1$ and a positive increasing function $Q:\rplus\rightarrow \rplus$ such that
\begin{eqnarray*} \norme{m}{U(t)}^2  & \leq & C\norme{m}{U(0)}^2 + t Q(\trinorme{m,t}{U}^2+\trinorme{m-1,t}{\derp{z}u\tang}^2+\trinorme{m,t}{F}^2+\trinorme{m-1,t}{\nabla F}^2) \\
 & & \hspace{15pt} +C \trinorme{\infty,t}{\derp{z}U}^2\int_0^t \norme{m-1}{\derp{z}\rho(s)}^2 + \norme{m-1}{\derp{z}u_3(s)}^2 ~ds. \end{eqnarray*}
Thus, under the conditions of \assuref{assu1}, we have, for $t\leq T^*$,
$$ \trinorme{m,t}{U}^2 + \eps \int_0^t \norme{m}{\nabla u}^2 \leq M_0 + Q(M+M_F)\left(t + \trinorme{\infty,t}{\derp{z}U}^2\int_0^t \norme{m-1}{\derp{z}\rho(s)}^2+\norme{m-1}{\derp{z}u_3(s)}^2~ds\right). $$ \label{hmeeinf} \end{esti}

\section[Proof of \thref{thEE}, part II: \textit{a priori} estimates on $\derp{z}u\tang$]{Proof of \thref{thEE}, part II \\ \textit{A priori} estimates on $\derp{z}u\tang$
 \sectionmark{Proof of \thref{thEE} part II: a priori estimates on $\derp{z}u\tang$}}
\sectionmark{Proof of \thref{thEE} part II: a priori estimates on $\derp{z}u\tang$}

\subsection{Conormal energy estimates}

In this section, we get estimates on the tangential components of $\derp{z}u$. We will perform conormal energy estimates on the first two components of the equation of the vorticity
$$ \omega=\rot u = \left( \begin{array}{c} \derp{y_{2}}u_{3}-\derp{z}u_{2} \\ \derp{z}u_{1}-\derp{y_{1}}u_{3} \\ \derp{y_{1}}u_{2}-\derp{y_{2}}u_{1} \end{array} \right). $$
By applying derivations to equation (\ref{NS}), we get that $w=\omega\tang$ solves the following equation:
\begin{equation} \rho\derp{t}w + \rho u\cdot\nabla w - \eps\mu\Delta w = {\cal M} , \label{roteqn} \end{equation}
$$ \mathrm{where} \hspace{20pt} {\cal M}=-\rho\omega\cdot\nabla u + \rho(\div u)\omega + \rot (\rho F) + {\cal M}^I + \eps {\cal M}^{II}. $$
The remainder ${\cal M}^I$ is a sum of terms written as $\derp{x_i}\rho(\derp{t}u_j+u\cdot\nabla u_j)e_j$, and ${\cal M}^{II}$ is a sum of terms written as $\kappa'(\rho) \derp{x_j}\rho (\derp{x_k,x_l}^2 u_j)e_j$, with $i,j,k,l\in\{1,2,3\}$, $i\neq j$, $(e_1,e_2,e_3)$ is the canonical basis of $\mathbb{R}^3$, and $\kappa$ is either $\lambda$ or $\mu$. 
The terms in ${\cal M}^{II}$ come from the derivation of the laplacian terms, but also from $\rot \sigma$, in which there are no terms with two derivatives on $\rho$ thanks to $\rot\nabla=0$. 
But the boundary condition on $\omega\tang$, according to (\ref{NBC}) and the non-penetration condition $u_{3}|_{z=0}=0$, is
$$ \mu(\rho)\omega\tang|_{z=0} = 2a u\tang^{\bot}|_{z=0}, $$
where for $v=(v_{1},v_{2})\in\rtwo$, $v^{\bot} = (-v_{2},v_{1})$, which makes integrations by parts difficult. So we introduce a modified vorticity:
$$ W = \omega\tang - \frac{2a}{\mu(\rho)} u\tang^{\bot} $$
By the tame estimate, we can harmlessly identify conormal Sobolev norms of $W$ with those of $\omega\tang$ and $\derp{z}u\tang$. The modified vorticity satisfies $W=0$ on the boundary, and solves the equation
\begin{equation} \rho\derp{t}W + \rho u\cdot\nabla W - \eps\mu\Delta W = H, \label{vorteqn} \end{equation}
$$ \mathrm{with} \hspace{10pt} H=2a[(\rho\derp{t}+\rho u\cdot\nabla - \eps\mu\Delta),\mu^{-1}] u\tang^\bot-2a\eps\mu^{-1}\lambda\nabla^{\bot}\tang\div u \hspace{70pt}  $$
$$ \hspace{70pt} + 2ak\mu^{-1}\rho^{\gamma-1}\nabla^{\bot}\tang\rho -2a\mu^{-1}\rho F\tang^\bot + {\cal M}\tang . $$
We now prove the following.

\begin{esti} The modified vorticity $W$ satisfies the following conormal energy estimate for every $m\geq 1$:
$$ c\left[\norme{m-1}{W(t)}^{2} + \eps\int_{0}^{t} \norme{m-1}{\nabla W(s)}^{2}~ds\right] \leq C\left[\norme{m-1}{W(0)}^{2} + \norme{m}{U(0)}^{2}\right] \hspace{110pt} $$
\begin{equation} \hspace{20pt} + Q(\trinorme{\Lip,t}{U}^2+\trinorme{1,\infty,t}{\nabla U}^2+\trinorme{\infty,t}{\nabla F}^2) \int_0^t (\norme{m}{U}^{2} + \norme{m-1}{\derp{z}U}^{2} + \norme{m}{F}^2 + \norme{m-1}{\nabla F\tang}^2)~ds, \label{hmvfinal} \end{equation}
for some positive increasing function $Q(z)$. \label{hmvort} \end{esti}

\preu{Proof:} we repeat the reasoning of section 4.1 on the above equation. We start by taking the $L^{2}$ scalar product of (\ref{vorteqn}) with $W$ and integrate by parts to get
$$ c\norme{0}{W(t)}^{2} + \eps\mu\int_{0}^{t}\norme{0}{\nabla W(s)}^{2}~ds \leq C\norme{0}{W(0)}^{2} \hspace{125pt} $$
\begin{equation} \hspace{80pt} +C\int_{0}^{t} \left[\trinorme{\infty,t}{\div(\rho u)}\norme{0}{W(s)}^{2} + \int_\Omega H(s)\cdot W(s) \right]~ ds. \label{l2vort} \end{equation}
We notice that, thanks to the compressibility equation (\ref{NSdiv}), $\trinorme{\infty,t}{\div(\rho u)}$ can be replaced by $\trinorme{\infty,t}{\derp{t}\rho}\leq \trinorme{1,\infty,t}{\rho}$. It remains to estimate $\int_\Omega H\cdot W$. 
The only terms we need this scalar product form for are the order-two terms in ${\cal M}^{II}$, and $\eps \int_\Omega \mu^{-1} (\Delta \mu) u\tang\cdot W$, which comes from the commutator $\eps \mu [\Delta,\mu^{-1}]u\tang$. After integrating $\int_\Omega \mu^{-1} (\Delta \mu) u\tang\cdot W$ by parts, we can simply bound $\mu$ and $\nabla \mu$ in $L^\infty$ and use Young's inequality on the terms involving $u$ and $W$. 
The term $\int_\Omega {\cal M}^{II} \cdot W$ is bounded by integration by parts on the variable $x_k$ with $k\neq 3$ and the use of Young's inequality with an adequate to absorb $\norme{0}{\nabla W}$ whenever it appears, while the terms containing $\derp{x_k,x_i}^2\rho$ are controlled by bounding the order-two part on $\rho$ in $L^\infty$ (thus getting $\
norme{1,\infty}{\nabla \rho}$). When $(k,l)=(3,3)$, we use (\ref{NS}) to replace $\eps\derp{zz}u$. For example,
$$ \eps\mu\derp{zz}u_1 = \rho\derp{t}u_1 + (\rho u \cdot\nabla) u_1 + \derp{y_1}P(\rho) - \rho F_1 - \eps [\mu (\derp{y_1 y_1}+\derp{y_2 y_2})u_1 + \lambda \derp{y_1}\div u] , $$
and the only difficulty is to control $\eps\int_0^t \norme{1}{\derp{z}u_1}^2$, but this is given by \propref{hmee}.

We can use the Cauchy-Schwarz inequality and handle the norm of the rest of $H$ as follows:
\begin{itemize}
\item the other commutators with $\mu^{-1}$ are controlled by $Q(1+\norme{\Lip}{U})(\norme{1}{U}+\norme{0}{W})$,
\item $\norme{0}{\rho^{\gamma-1}\nabla\tang^{\bot}\rho} \leq C\norme{1}{U}$ as $\rho$ is assumed to be uniformly bounded,
\item $\norme{0}{\rho\omega\cdot\nabla u\tang + \rho(\div u)\omega} \leq C\norme{\infty}{\nabla u}(\norme{0}{W}+\norme{0}{u})$,
\item ${\cal M}^I$ is bounded using $\norme{0}{\derp{x_{i}}\rho(\derp{t}u_j+u\cdot\nabla u_j)} \leq \trinorme{\infty,t}{\nabla \rho}(\norme{1}{u}+\trinorme{\infty,t}{\nabla u}\norme{0}{u})$,
\item and finally, we have $\eps\norme{0}{\nabla\tang \div u}\norme{0}{W} \leq \eps(\norme{1}{\nabla u}^{2}+\norme{0}{W}^{2})$ and use \estiref{hmee} to control $\eps\norme{1}{\nabla u}$,
$$ \eps\int_{0}^{t}\norme{1}{\nabla u(s)}^{2}~ds \leq C\left(\norme{1}{U(0)}^{2}+(1+\trinorme{\Lip,t}{U}^{2})\int_{0}^{t}(\norme{1}{U(s)}^{2} + \norme{0}{\nabla U(s)}^{2})~ds\right). $$
\end{itemize}

We thus get the final $L^{2}$ estimate on $W$:
$$ c\left[\norme{0}{W(t)}^{2}+\eps\int_{0}^{t}\norme{0}{\nabla W(s)}^{2}~ds\right] \leq C(\norme{0}{W(0)}^{2} + \norme{1}{U(0)}^{2}) \hspace{120pt} $$
$$ \hspace{60pt} +Q(1+\trinorme{\Lip,t}{U}^{2}+\trinorme{1,\infty,t}{\nabla U}^2)\int_{0}^{t} \norme{1}{U}^{2} + \norme{0}{\derp{z}U}^{2}+\norme{1}{F}^2 + \norme{0}{\nabla F}^2~ds . $$

Now we move on to the $H^{m-1}_{co}$ estimate with $m>1$, which will also follow the same pattern as the previous section. We apply $Z^{\alpha}$ to (\ref{vorteqn}), for $|\alpha|\leq m-1$ and isolate the maximum order terms:
$$ \rho\derp{t}Z^{\alpha}W + \rho u\cdot \nabla Z^{\alpha} W - \eps\mu\Delta Z^{\alpha}W = Z^{\alpha}H - {\cal C}^{\alpha}_{W}, $$
where ${\cal C}^{\alpha}_{W}= [Z^{\alpha},\rho]\derp{t}W+[Z^{\alpha},\rho u\cdot\nabla]W - \eps[Z^\alpha,\mu\Delta]W$ contains the commutators. Multiplying the above equation by $Z^\alpha W$ in $L^{2}$, we have
$$ c\left[\norme{0}{Z^{\alpha}W(t)}^{2} + \eps\int_{0}^{t}\norme{0}{Z^{\alpha}\nabla W}^{2}~ds\right] \leq C\left[\norme{0}{Z^\alpha W(0)}^{2} + \trinorme{1,\infty,t}{\rho}\int_0^t \norme{m-1}{W}^2~ds\right] \hspace{50pt} $$
\begin{equation} \hspace{60pt} + C\eps\int_{0}^{t} \norme{m-2}{\nabla W}^2~ds + C\int_0^t \left[\norme{m-1}{H}^2+\int_\Omega |{\cal C}^{\alpha}_{W}\cdot Z^\alpha W|\right] ~ds . \label{hmvort1} \end{equation}

The terms on the second line of this inequality are the ones we need to control. 
First, $\eps\int_0^t \norme{m-2}{\nabla W}^2$, which comes from changing $\nabla Z^\alpha W$ into $Z^\alpha \nabla W$ (using \propref{comprop}) in the integration by parts on the laplacian, is dealt with by induction, by using (\ref{hmvfinal}) at a lower order. 
Then, $\norme{m-1}{H}$ is easily estimated as above, also using the tame estimate in \propref{tame}:
\begin{eqnarray*} \int_0^t \norme{m-1}{H(s)}^2~ds & \leq & \eps\int_0^t \norme{m}{\nabla u}^{2}~ds + Q(\trinorme{\infty,t}{U}^2+ \trinorme{1,\infty,t}{\nabla U}^2+\trinorme{\infty,t}{\nabla F}^2) \\
 & & \hspace{10pt} \times \left(\int_0^t \norme{m}{U}^{2}+\norme{m-1}{\derp{z}U}^{2} + \norme{m}{F}^2+ \norme{m-1}{\nabla F\tang}^2~ds\right) . \end{eqnarray*}
Using \estiref{hmee} to cover the worst term $\eps\int_{0}^{t}\norme{m}{\nabla u(s)}^{2}~ds$, we have
\begin{eqnarray} \int_0^t \norme{m-1}{H(s)}^2~ds & \leq & C\norme{m}{U(0)}^2 + Q(\trinorme{\infty,t}{U}^2+\trinorme{1,\infty,t}{\nabla U}^2+\trinorme{\infty,t}{\nabla F}^2) \nonumber \\
 & & \hspace{10pt} \times \int_0^t \norme{m}{U}^2+\norme{m-1}{\derp{z}U}^2 + \norme{m}{F}^2 + \norme{m-1}{\nabla F}^2~ds . \label{hmvrhs} \end{eqnarray}

We finally need to estimate the commutators in $\norme{0}{{\cal C}^{\alpha}_{W}}$. The first and last terms of ${\cal C}^{\alpha}_{W}$ are easily estimated as in the previous section,
\begin{equation} \int_0^t \norme{0}{[Z^{\alpha},\rho]\derp{t}W}^2~ds \leq C\int_0^t (\trinorme{\infty,t}{\rho}^2+\trinorme{\infty,t}{\derp{t}W}^2)(\norme{m-1}{W}^2+\norme{m-1}{U}^2)~ds \label{hmcom1} \end{equation}
by \propref{funcprop}, and again using the decomposition in \propref{comprop} (b) and integrating by parts, we have
\begin{equation} \eps \int_{\Omega} |[Z^\alpha,\mu\Delta]W\cdot Z^{\alpha} W| \leq \frac{1+\trinorme{\Lip,t}{\rho}}{c}\eps\norme{m-2}{\derp{z}W}^{2}+\frac{c}{4}\eps\norme{m-1}{\derp{z}W}^{2}, \label{hmcom2} \end{equation}
in which the first term is estimated by induction (using (\ref{hmvfinal}) at order $m-2$) and the second is absorbed by the left-hand side of (\ref{hmvort1}). 
Finally, instead of applying the commutator results we have used so far to $\int_0^t \norme{0}{[Z^\alpha, \rho u\cdot\nabla] W}^2$, which would yield a $\norme{m-2}{\derp{zz}^2 U}$ term that we do not expect to control uniformly in $\eps$, we use \propref{comlem}:
$$ \int_{0}^{t}\norme{0}{[Z^{\alpha},(\rho u)(s)\cdot\nabla]W(s)}^{2}~ds \leq \int_{0}^{t}Q(\trinorme{\Lip,t}{U}^{2}+\trinorme{1,\infty,t}{\nabla U}^{2})\norme{m-1}{W(s)}^{2}~ds $$
$$ \hspace{135pt} +\int_{0}^{t}Q(\trinorme{1,\infty,t}{\nabla U}^{2})\norme{m-1}{\nabla U(s)}^{2}~ds. $$
This finishes off the proof of \estiref{hmvort}. $\square$

\subsection{$L^\infty$ estimates}

To deal with $\norme{1,\infty}{\derp{z}u\tang}$, we examine $\omega\tang$. The result is:
\begin{esti} There exists a positive, increasing function on $\rplus$, $Q$, such that
$$ \norme{1,\infty}{\omega\tang(t)}^2 \leq Q(\norme{6}{U(0)}^2+\norme{5}{\omega\tang(0)}^2+\norme{1,\infty}{\omega\tang(0)}^2) \hspace{170pt} $$
$$ \hspace{30pt} + Q(\trinorme{\Lip,t}{U}^2+\trinorme{1,\infty,t}{\nabla U}^2+{\cal N}_6(t,F))\int_0^t 1+\norme{6}{U(s)}^2+\norme{6}{\omega\tang(s)}^2+\norme{5}{\derp{z}\rho(s)}^2~ds. $$  \label{linfvort} \end{esti}

We can now deduce the following update of \propref{hmeeinf}.
\begin{coro} Under the conditions of \assuref{assu1}, for $m\geq 7$ and $t\leq T^*$,
$$ \trinorme{m,t}{U}^2+\trinorme{m-1,t}{\nabla u\tang}^2+\trinorme{1,\infty,t}{\nabla u\tang}^2 + \eps\int_0^t \norme{m}{\nabla u}^2 + \norme{m-1}{\nabla W}^2~ds \hspace{70pt} $$
$$ \hspace{80pt} \leq Q(M_0) + Q(M+M_F)\left(t+ \trinorme{1,\infty,t}{\derp{z}(\rho,u_3)}^2 \int_0^t \norme{m-1}{\derp{z}\rho}^2+\norme{m-1}{\derp{z}u_3}^2~ds\right). $$
\label{uvortbd} \end{coro}

The first tool to prove \propref{linfvort} is the following version of the maximum principle.
\begin{propo} Consider $X$ a Lipschitz-class solution to the following hyperbolic-parabolic system on the half-space $\Omega\subset\mathbb{R}^3$:
\begin{equation} \left\{ \begin{array}{rcll} a\derp{t}X + b\cdot\nabla X - \eps\mu(a)\Delta X & = & G & \mathrm{in}~\Omega \\
X & = & h & \mathrm{on}~\dOmega \\ X|_{t=0} & = & X_0\in L^\infty(\Omega), & \end{array} \right. \label{mpsys} \end{equation}
in which $X:\rplus\times\Omega\rightarrow \mathbb{R}^d$, $a$ is a scalar function, bounded and positive uniformly in time ($\inf_{(t,x)} a(t,x)>c>0$), $b:\rplus\times\Omega\rightarrow\mathbb{R}^3$ is tangent to the boundary, $G:\rplus\times\Omega\rightarrow\mathbb{R}^d$ and $h:\rplus\times \rtwo \rightarrow \mathbb{R}^d$ are in $L^\infty$ and $\mu$ is a positive regular function of $a$. 
We assume that $a$ and $b$ satisfy $\derp{t}a+\div b = 0$. Then we have, for $t\in [0,T]$,
$$ \norme{\infty}{X(t)} \leq \norme{\infty}{X(0)} + \trinorme{\infty,T}{h} + \frac{1}{\inf a}\int_0^t \norme{\infty}{G(s)} + \eps\norme{\infty}{\mu'(a(s))\nabla a(s) \cdot \nabla X(s)} ~ds. $$ \label{mp1} \end{propo}

Note that when $\mu$ is constant, the term involving $\nabla(\mu(a))$ vanishes on the right-hand side. 
In fact, in that case we can directly apply the result to the modified vorticity $W$ instead of $\omega$, which has the further advantage of satisfying a homogeneous boundary condition ($h=0$). 
In the non-constant case however, if we consider $W$, we cannot deal with the term $\derp{zz}\rho$ that emerges from the commutator $[\Delta,\mu^{-1}]$, so we will use the proposition on $\omega\tang$. 
Having $\nabla X$ in the right-hand side looks disastrous, but in the context of the estimation process, the factor $\eps$ intervenes crucially.
\newline

\preu{Proof:} we look for a function $g(t)$ that will control $\norme{\infty}{X(t)}$. To do so, we perform energy estimates on the functions $(X-g)_+$ and $(X+g)_-$ ($g$ will indifferently designate the scalar function $g$ and the vector $g\mathbf{1}\in\mathbb{R}^d$), with the convention that, for a given scalar function $f$, $f=f_+ + f_-$ (so $f_-\leq 0$).

We concentrate on $(X-g)_+$, the procedure on $(X+g)_-$ being identical. We want to find a function $g(t)$ that satisfies
$$ \norme{0}{(X(0)-g(0))_+} = 0 \mathand \frac{d}{dt}\norme{0}{(X(t)-g(t))_+}^2 \leq 0,~t> 0; $$
such a function $g$ is an upper bound for $\norme{\infty}{X_+}$. Choosing $g(0)=\norme{\infty}{X(0)}$, we examine
\begin{eqnarray} \frac{1}{2}\frac{d}{dt}\left(\int_\Omega a|(X(t)-g(t))_+|^2\right) & = & \int_\Omega a\derp{t}X\cdot(X-g)_+ - \int_\Omega a\derp{t}g\cdot(X-g)_+ \nonumber \\
 & & \hspace{10pt} + \frac{1}{2} \int_\Omega \derp{t}a |(X-g)_+|^2 \nonumber \\
 & = & \int_\Omega (-b\cdot\nabla X + \eps\mu\Delta X + G)\cdot(X-g)_+ \nonumber \\
 & & \hspace{10pt} + \int_\Omega -a\derp{t}g\cdot(X-g)_+ + \frac{1}{2}\derp{t}a|(X-g)_+|^2 \label{mpee} \end{eqnarray}
Because of the scalar product with $(X-g)_+$, which is zero wherever $X-g\leq 0$, we can replace $\nabla X$ by $\nabla ((X-g)_+)$ and $\Delta X$ by $\Delta((X-g)_+)$. Then, integration by parts gives us
$$ -\int_\Omega b\cdot\nabla X\cdot(X-g)_+ = \frac{1}{2}\int_\Omega \div b |(X-g)_+|^2 = -\frac{1}{2} \int_\Omega \derp{t}a|(X-g)_+|^2, $$
by the first line of (\ref{mpsys}), so this term cancels out with the final term of (\ref{mpee}), and it remains to look at the term involving the laplacian. When integrating by parts, we need to guarantee that $(X-g)_+$ vanishes on the boundary. We therefore impose the condition
\begin{equation} g(t) \geq \trinorme{L^\infty(\dOmega),T}{h} \label{gbc} \end{equation}
and write
\begin{eqnarray*} \eps\mu\int_\Omega \Delta X\cdot(X-g)_+ & = & -\eps\mu\norme{0}{\nabla ((X-g)_+)}^2 - \eps \int_\Omega \nabla (\mu(a))\cdot\nabla X\cdot (X-g)_+ \\
 & \leq & -\eps \int_\Omega \mu'(a) \nabla a\cdot \nabla X\cdot (X-g)_+ \end{eqnarray*}
Therefore, setting $\til{G} = G-\eps \mu'(a) \nabla X\cdot\nabla a$, of (\ref{mpee}) there only remains
$$ \frac{d}{dt}\left(\int_\Omega a|(X(t)-g(t))_+|^2\right) \leq 2\int_\Omega (\til{G}-a\derp{t}g)\cdot(X-g)_+, $$
which is negative if, for each $j\in\{1,\cdots,p\}$, $\til{G}_j-a\derp{t}g\leq 0$ on $\Omega$. Integrating this leads to $g(t)\geq g(0)+\int_0^t \frac{\til{G}_j(s,x)}{a(s,x)}~ds$ for every $x\in\Omega$.

Identical estimates on $(X+g)_-$ lead to $g(t)\geq g(0)-\int_0^t \frac{\til{G}_j(s,x)}{a(s,x)}~ds$, so, also taking into account (\ref{gbc}), we choose
$$ g(t) = \norme{\infty}{X(0)} + \trinorme{L^{\infty}(\dOmega),T}{h} + \int_0^t \norme{\infty}{\frac{\til{G}(s)}{a(s)}}~ds, $$
and this controls $\norme{\infty}{X(t)}$ as desired. $\square$
\newline

\preu{Proof of \propref{linfvort}:} we immediately apply \propref{mp1} to $w:=\omega\tang$ and $Z_j w$, to get
$$ \norme{1,\infty}{w(t)} \leq \norme{1,\infty}{w(0)} + 2a\trinorme{1,\infty,t}{\mu^{-1} u\tang} + \int_0^t \norme{1,\infty}{\frac{{\cal M}\tang(s)}{\rho(s)}} + \sum_{j=0}^3\norme{\infty}{{\cal C}^j_w(s)}~ds, $$
where ${\cal C}_w^j=[Z_j,\rho]\derp{t}w+[Z_j,\rho u\cdot\nabla] w - \eps[Z_j,\mu\Delta]w$ are the commutator terms that appear in the equations on $Z_j w$. Thus, by using the Cauchy-Schwarz and Young inequalities,
\begin{equation} \norme{1,\infty}{w(t)}^2 \leq \norme{1,\infty}{w(0)}^2 + 4a^2\trinorme{1,\infty,t}{\mu u\tang}^2 + t\int_0^t \norme{1,\infty}{\frac{{\cal M}\tang(s)}{\rho(s)}}^2+\sum_{j=0}^3 \norme{\infty}{{\cal C}^j_w(s)}^2~ds . \label{inftyv1} \end{equation}
Note that a factor $t$ has been extracted in front of the source terms and commutators, so we can be satisfied with fairly crude bounds on these. The $W^{1,\infty}_{co}$ term stemming from the boundary is easily dealt with using \propref{sobo}.
\newline

Let us start with the control of the commutators. Given that, by (\ref{commul}),
$$ \eps\mu[Z_3,\derp{zz}]w=-2\eps\mu\phi'\derp{zz}w-\eps\mu\phi''\derp{z}w , $$
we replace $\eps\mu\derp{zz}w$ by its expression in the equation. So we write, using the notation $\Delta\tang = \derp{y_1 y_1}+\derp{y_2 y_2}$,
\begin{eqnarray*} \eps^2 \mu^2\norme{\infty}{[Z_3,\derp{zz}]w}^2 & \leq & C\eps^2\norme{\infty}{\derp{z}w}^2 + C\norme{\infty}{\rho\derp{t}w+\rho u\cdot\nabla w - \eps\mu\Delta\tang w - {\cal M}\tang}^2 \\
 & \leq & C\eps^2\norme{\infty}{\derp{z}w}^2 + C(\norme{1,\infty}{\rho}^2+\norme{\Lip}{u}^2\norme{1,\infty}{w}^2) \\
 & & \hspace{30pt} + C\eps^2\norme{2,\infty}{w}+C\norme{\infty}{{\cal M}\tang}^2, \end{eqnarray*}
by using the fact that $\norme{\infty}{u_3\derp{z}w} \leq C\norme{\infty}{\derp{z}u_3}\norme{\infty}{Z_3 w}$. 
We see that $\norme{\infty}{{\cal C}_w^3}$ is bounded, among other terms, by $\norme{\infty}{{\cal M}\tang}$; we'll examine the source term last. 
In this commutator, it remains to control
$$ J:=\int_0^t \eps^2 (\norme{\infty}{\derp{z}w(s)}^2+\norme{2,\infty}{w(s)}^2)~ds . $$
On both parts of this term, we start by applying \propref{sobo}, the Sobolev embedding theorem, to get
$$ J \leq C\int_0^t \eps^{2}(\norme{3}{\derp{zz}w(s)}^2+\norme{5}{\derp{z}w(s)}^2+\norme{5}{w(s)}^2)~ds. $$
The second term, $\int_0^t \eps^2\norme{5}{\derp{z}w}^2$, is controlled by using \estiref{hmvort}, and, like above, the first term, $\int_0^t \norme{3}{\eps\derp{zz}w}^2$, is controlled by replacing $\eps\derp{zz}w$ by its expression in (\ref{roteqn}), and by using the tame estimate:
\begin{eqnarray*} \int_0^t \norme{3}{\eps\derp{zz}w(s)}^2~ds & \leq & C \int_0^t (\trinorme{\infty,t}{\rho}^2+\trinorme{\infty,t}{\rho u}^2)\norme{4}{w}^2+\trinorme{1,\infty,t}{w}^2\norme{3}{U}^2 \\
& & \hspace{50pt} + \norme{3}{\rho u_3 \derp{z} w} + \norme{5}{w}^2 + \norme{3}{{\cal M}\tang}^2 ~ ds \end{eqnarray*}
The source term norm $\int_0^t \norme{3}{{\cal M}\tang(s)}^2~ds$ is given by (\ref{hmvrhs}), and $\int_0^t \norme{3}{\rho u_3 \derp{z} w}^2~ds$ is controlled in the same way as is done in \propref{comlem}:
\begin{eqnarray*} \int_0^t \norme{3}{(\rho u_3 \derp{z} w)(s)}^2~ds & = & \int_0^t \norme{3}{\frac{(\rho u_3)(s)}{\phi} Z_3 w(s)}^2~ds \\
 & \leq & C \int_0^t \trinorme{1,\infty,t}{w}^2 \norme{4}{\derp{z}(\rho u)}^2+(\trinorme{\Lip,t}{\rho u}^2+\trinorme{1,\infty,t}{\nabla (\rho u)}^2)\norme{4}{w}^2~ds \end{eqnarray*}
The commutator $\eps[Z_3,\mu]\Delta w=Z_3(\mu(\rho)) \Delta w$ is controlled by replacing $\eps\derp{zz}w$ by its expression. Therefore, in total, we have
$$ \eps^2\int_0^t \norme{\infty}{[Z_3,\mu\Delta]w}^2~ds \leq C[\norme{6}{U(0)}^2+\norme{5}{w(0)}^2] + C\int_0^t \norme{\infty}{{\cal M}\tang(s)}^2~ds \hspace{50pt} $$
\begin{equation} \hspace{60pt} + C\int_0^t Q(\trinorme{\Lip,t}{U}^2+\trinorme{1,\infty,t}{\nabla U}^2)(1+\norme{6}{U}^2+\norme{5}{w}^2+\norme{5}{\derp{z}\rho}^2)~ds+{\cal N}_6(t,F) \label{lapcomwinf} \end{equation}
It is now straightforward to notice that the whole commutator satisfies (\ref{lapcomwinf}), the main tool to bound $\norme{\infty}{[Z_j,\rho]\derp{t}w+[Z_j,\rho u \cdot\nabla]w}$ being property (\ref{hardyinf}), as used in \propref{comlem}.
\newline

Moving on to the source term, we now focus on the terms in $\eps {\cal M}^{II}$ which arise when the viscosity coefficients are not constant. First, we examine terms involving one normal derivative on $u_1$ or $u_2$. These are linked to $\nabla w$, which we split into two parts:
\begin{eqnarray*} \eps \nabla w & = & \eps\nabla ((\rot u)_\tau) = \eps \nabla (\derp{z} u_\tau^\bot - \nabla_\tau^\bot u_3) \\
 & = & \eps \left(\begin{array}{c} \nabla_\tau \\ \derp{z} \end{array} \right) \derp{z}u_\tau^\bot - \eps \nabla (\nabla_\tau^\bot u_3) \end{eqnarray*}
We are interested in the $L^2$-in-time $W^{1,\infty}_{co}$ norm of this. On the terms with only one normal derivative, which are the second and the $\nabla_\tau$ components of the first, we apply the Sobolev embedding, \propref{sobo}, leading to having to bound $\norme{5}{\eps \derp{zz}u_j}^2$ for a certain $j$; once again we can do so by replacing $\derp{zz}u_j$ by its expression in (\ref{NS}). 
The remains of the first term, $\derp{zz}u_\tau$, are no more difficult: having replaced $\eps\derp{zz}u_\tau$ using (\ref{NS}), all the subsequent terms are dealt with easily using the Sobolev embedding in most places, including on $\eps\lambda \nabla_\tau \derp{z}u_3$, which comes from the $\nabla\div$term, and another replacement of $\eps\derp{zz}u_3$ completes the estimate.

In the remaining terms of $\eps {\cal M}^{II}$, terms with two normal derivatives on $u$ can appear: we then replace them using the equation. Terms with two conormal derivatives are dealt with by using the Sobolev embedding inequality.
\newline

Estimating the $L^2$-in-time conormal-Lipschitz norm of the rest of the source term, ${\cal M}\tang$, which, we remind the reader, is the tangential component of $-\rho\omega\cdot\nabla u + \rho(\div u)\omega + \rot (\rho F) + {\cal M}^I$, is straight-forward. $\square$

\section[Proof of \thref{thEE}, part III: \textit{a priori} estimates on $\derp{z}u_3$]{Proof of \thref{thEE}, part III \\ \textit{A priori} estimates on $\derp{z}u_3$
 \sectionmark{Proof of \thref{thEE} part III: a priori estimates on $\derp{z}u_3$}}
\sectionmark{Proof of \thref{thEE} part III: a priori estimates on $\derp{z}u_3$}

To begin with, we bound $\int_0^t \norme{m-1}{\derp{z}u_3(s)}^2~ds$, which will essentially only need the tame estimate in \propref{tame}. The compressibility equation (\ref{NSdiv}) gives us
$$ \rho Z^\alpha \derp{z} u_3 = - Z^\alpha \derp{t}\rho - Z^\alpha(u\cdot\nabla \rho) - Z^\alpha(\rho \divt u\tang) - [Z^\alpha,\rho]\derp{z}u_3, $$
with $|\alpha|\leq m-1$ and introducing the notation $\divt u\tang = \derp{y_1}u_1 + \derp{y_2}u_2$. Integrating the square of this equality in time and space, we quickly get
$$ c\int_0^t \norme{0}{Z^\alpha \derp{z}u_3(s)}^2~ds \leq C (1+\trinorme{\Lip,t}{U}^2) \int_0^t \norme{m}{U(s)}^2~ds \hspace{75pt} $$
$$ \hspace{50pt} + C\trinorme{\infty,t}{u_3}^2 \int_0^t \norme{m-1}{\derp{z}\rho(s)}^2~ds + C\trinorme{1,\infty,t}{\rho}^2 \int_0^t \norme{m-2}{\derp{z}u_3(s)}^2~ds, $$
using \propref{funcprop} on the commutator term. With this, we can perform an induction on $m>0$ and conclude that there exists a positive increasing polynomial function $Q$ such that
$$ \int_0^t \norme{m-1}{\derp{z}u_{3}(s)}^2~ds \leq Q(1+\trinorme{\Lip,t}{U}^2) \int_0^t \norme{m}{U(s)}^2 + \norme{m-1}{\derp{z}\rho(s)}^2~ds . $$
Essentially, this means that $H^{m-1}_{co}([0,T]\times\Omega)$ norms of $\derp{z}u_3$ can be replaced by the same norms of $\derp{z}\rho$ and $W^{1,\infty}_{co}$ norms of $\derp{z}u_3$, which, in turn, will be bounded by $H^{m_0}_{co}$ norms for a certain $m_0$, so we need to extract a small parameter, in this case $t$, to close the estimate. We already have
\begin{equation} \int_0^t \norme{m-1}{\derp{z}u_3(s)}^2~ds \leq Q(1+\trinorme{\Lip,t}{U}^2)\left(t\trinorme{m,t}{U} + \int_0^t \norme{m-1}{\derp{z}\rho(s)}^2~ds \right) , \label{hmdzu3} \end{equation}
so we need the estimates on $\norme{1,\infty}{\derp{z}u_3}$ and $\norme{1,\infty}{\derp{z}\rho}$ to yield a small parameter, $t$ or $\eps$, as we cannot use (\ref{l2time}) on the $L^2$-in-time norm present here (by \estiref{linfu}, $\norme{1,\infty}{U}$ already yields the desired parameter).
\newline

For now, we focus on $\norme{1,\infty}{\derp{z}u_3}$. Simply reading (\ref{NSdiv}), we have
$$ \trinorme{\infty,t}{\derp{z}u_3}^2 \leq \trinorme{\infty,t}{\derp{t}\rho + \divt (\rho u\tang)}^2 + \trinorme{\infty,t}{u_3}^2\trinorme{\infty,t}{\derp{z}\rho}^2, $$
to which we can apply (\ref{l2timeinf}), and we get
\begin{eqnarray} \trinorme{\infty,t}{\derp{z}u_3}^2 & \leq & Q(1+\norme{\Lip}{U(0)}^2) + \int_0^t Q(1+\norme{2,\infty}{U(s)}^2+\norme{1,\infty}{\derp{z}\rho(s)}^2)~ds \nonumber \\
 & \leq & Q(1+\norme{\Lip}{U(0)}^2) + tQ(1+\trinorme{5,t}{U}^2+\trinorme{4,t}{\derp{z}U}^2 + \trinorme{1,\infty,t}{\derp{z}\rho}^2) \label{linfnd1} \end{eqnarray}
We can now do the same for $\norme{1,\infty}{\derp{z}u_3}$, bearing in mind that there is a commutator, $[Z_j,\rho]\derp{z}u_3 = (Z_j\rho)\derp{z}u_3$, and that we do not want to lose derivatives on $\norme{1,\infty}{\derp{z}\rho}$ this time:
\begin{eqnarray*} \trinorme{1,\infty,t}{\derp{z}u_3}^2 & \leq & \trinorme{1,\infty,t}{\derp{t}\rho + \divt(\rho u\tang )}^2 + \trinorme{1,\infty,t}{u_3}^2\trinorme{1,\infty,t}{\derp{z}\rho}^2 + \trinorme{1,\infty}{\rho}^2\trinorme{\infty}{\derp{z}u_3}^2 \\
 & \leq & Q(1+\norme{2,\infty}{U(0)}^2) + \int_0^t \norme{2,\infty}{\derp{t}\rho+ \divt (\rho u\tang )}^2 \\
 & & \hspace{30pt} +\trinorme{1,\infty,t}{\derp{z}\rho}^2\left(\norme{1,\infty}{u_3(0)}^2+\int_0^t \norme{2,\infty}{u_3(s)}~ds\right) \\
 & & \hspace{30pt} + \trinorme{\infty,t}{\derp{z}u_3}^2\left(\norme{1,\infty}{\rho(0)}^2+\int_0^t \norme{2,\infty}{\rho(s)}^2~ds\right) \\
 & \leq & Q(1+\norme{\Lip}{U(0)}^2+\norme{1,\infty}{\nabla U(0)}^2) + C\trinorme{1,\infty,t}{\derp{z}\rho}^4 \\
 & & \hspace{30pt}  + tQ(\trinorme{6,t}{U}^2+\trinorme{5,t}{\derp{z}U}^2+\trinorme{1,\infty,t}{\derp{z}\rho}^2) , \end{eqnarray*}
by (\ref{linfnd1}), which, by applying (\ref{l2time}) provides us with the following.

\begin{esti} (a) There exists an increasing, positive polynomial function $Q$ on $\rplus$ such that
$$ \trinorme{1,\infty,t}{\derp{z}u_3}^2 \leq Q(1+\norme{\Lip}{U(0)}^2+\norme{1,\infty}{\nabla U(0)}^2) + Q(\trinorme{1,\infty,t}{\derp{z}\rho}^2) \hspace{110pt} $$
\begin{equation} \hspace{40pt} + t Q\left( 1+\trinorme{6,t}{U}^2+\trinorme{5,t}{\nabla u\tang}^2+\trinorme{1,\infty,t}{\derp{z}\rho}^2+\int_0^t \norme{6}{\derp{z}u_3(s)}^2+\norme{6}{\derp{z}\rho(s)}^2~ds \right) \label{linfnd} \end{equation}

(b) Combining (\ref{hmdzu3}) and (\ref{linfnd}), we get that, for $m\geq 7$,
$$ \int_0^t \norme{m-1}{\derp{z}u_3(s)}^2~ds \leq Q(1+\norme{\Lip}{U(0)}^2+\norme{1,\infty}{\nabla U(0)}^2) \hspace{150pt} $$
$$ \hspace{35pt} + t Q\left(1+\trinorme{m,t}{U}^2+\trinorme{m-1,t}{\nabla u\tang}^2 + \trinorme{1,\infty,t}{\nabla u\tang}^2 + \int_0^t \norme{m-1}{\derp{z}u_3(s)}^2~ds \right) $$
\begin{equation} + Q\left(\trinorme{1,\infty,t}{\derp{z}\rho}^2 + \int_0^t \norme{m-1}{\derp{z}\rho(s)}^2~ds \right) \hspace{130pt} \label{hmnddef} \end{equation}
\label{hmnd} \end{esti}

As a result, we can update \cororef{uvortbd}. Set
\begin{eqnarray*} \til{E}^\eps _{m}(t,U)^2 & = & \trinorme{m}{U}^2+\trinorme{m-1,t}{\nabla u\tang}^2 + \int_0^t \norme{m-1}{\derp{z}u_3(s)}^2~ds+\trinorme{1,\infty,t}{\nabla u}^2 \\
 & & \hspace{40pt} + \eps\int_0^t \norme{m}{\nabla u}^2 + \norme{m-1}{\nabla W}^2~ds . \end{eqnarray*}
This contains ${\cal E}_m(t,U)$ except the terms involving $\derp{z}\rho$, plus the $W^{1,\infty}_{co}$ norm of $\derp{z}u_3$ and the gradient terms from the energy estimates. The combination of \cororef{uvortbd} and \propref{hmnd} (b) gives us the following.
\begin{coro} Under the conditions of \assuref{assu1}, for $m\geq 7$ and $t\leq T^*$,
$$ \til{E}^\eps_m(t,U) \leq Q(M_0) + Q(M+M_F)\left(t + Q\left(\trinorme{1,\infty,t}{\derp{z}\rho}^2 + \int_0^t \norme{m-1}{\derp{z}\rho(s)}^2~ds\right)\right) $$ \label{endp3} \end{coro}
We see on this estimate that it remains to look at the normal derivative of the density.

\section[Proof of \thref{thEE}, part IV: \textit{a priori} estimates on $\derp{z}\rho$]{Proof of \thref{thEE}, part IV \\ \textit{A priori} estimates on $\derp{z}\rho$
 \sectionmark{Proof of \thref{thEE} part IV: a priori estimates on $\derp{z}\rho$}}
\sectionmark{Proof of \thref{thEE} part IV: a priori estimates on $\derp{z}\rho$}

\subsection{Conormal energy estimates}

In this section, we examine $R:=\derp{z}\rho$. The equation satisfied by $R$ is given by the differentiation of (\ref{NSdiv}):
$$ \derp{t}R + u\cdot\nabla R + R(\div u+\derp{z}u_3) + \rho\derp{zz}u_3 = - \rho\derp{z}\divt u\tang - \derp{z}u\tang \cdot\nabla\tang \rho. $$
A very problematic term appears in this equation: $\derp{zz}u_3$. The idea is to multiply the equation by $l\eps:=(\mu+\lambda)\eps$, in order to replace $l\eps\derp{zz}u_3$ by its expression in the equation,
\begin{equation} l\eps\derp{zz} u_3 = \rho \derp{t}u_3 + \rho u\cdot\nabla u_3 - \mu\eps\Delta\tang u_3 - \lambda\eps\derp{z}\divt u\tang + R. \label{NSrwR} \end{equation}
This brings us to
\begin{equation} l\eps(\derp{t}R + u\cdot\nabla R) + \rho P'(\rho) R = h, \label{eqnR} \end{equation}
where we remind the reader that $\rho P'(\rho) = k\rho^\gamma$ (and $P'$ is positive; we can thus expect our results to extend to other, strictly increasing and positive pressure laws), and $h$ will be treated as a source term:
\begin{eqnarray*} h & = & \eps[\rho(\mu\Delta\tang u_3+\lambda\derp{z}\divt u\tang+\sigma(\nabla U)_3) \\
 & & \hspace{20pt} -l\derp{z}u\tang\cdot\nabla\tang\rho-l(\div u+\derp{z}u_3)R]-\rho^2(\derp{t}u_3+u\cdot\nabla u_3)+\rho^2 F_3. \end{eqnarray*}

In these estimates, we will strongly use \assuref{assu1}: for $t\leq T^*$, we have $\rho(t,x)\geq c_0>0$, $|u_3(t,x)|$ is uniformly bounded near the boundary on the same time interval, and ${\cal E}_m(t,U)\leq M$ for some $m\geq 7$. We use the latter assumption to simplify the presentation of what follows.

\begin{esti} Under the conditions of \assuref{assu1}, for $m\geq 1$, there exists a positive, increasing function on $\rplus$, $Q$, such that, for $t\leq T^*$ and $\eps\leq \eps_0$,
$$ \int_0^t \norme{m-1}{R(s)}^2~ds \leq Q(\norme{m-1}{R(0)}^2 +\norme{m}{U(0)}^2) + (t+\eps)Q(M+M_F) + Q(\trinorme{1,\infty,t}{R}^2). $$
\label{hmeeR} \end{esti}

\preu{Proof:} we start with the $L^2$ estimate. We multiply (\ref{eqnR}) by $R$ and integrate in space, and, as usual, we integrate the term containing $u\cdot\nabla R$ by parts and get
\begin{equation} \frac{l\eps}{2}\frac{d}{dt}(\norme{0}{R}^2) + \gamma \int_\Omega \rho P'(\rho) R^2 = \int_\Omega \left(\frac{l\eps}{2} R^2 \div u + \frac{l'(\rho)\eps}{2} R^2 u\cdot\nabla \rho + hR\right). \label{l2eeR} \end{equation}
The coefficients $l$ and $\rho P'$ are bounded from below, so there exist $0<c<1<C$ such that the left-hand side, integrated in time, is greater than
$$ c\left(\eps \norme{0}{R(t)}^2 + \int_0^t \norme{0}{R(s)}^2~ds\right) -C\norme{0}{R(0)}^2 . $$
Then, using the Young inequality on the right-hand side of (\ref{l2eeR}) with the parameter $\eta=\frac{c}{2}$, we get
$$ c\eps \norme{0}{R(t)}^2 + \frac{c}{2}\int_0^t \norme{0}{R(s)}^2~ds \leq C\norme{0}{R(0)}^2 + C\int_0^t \eps \trinorme{\Lip,t}{U}\norme{0}{R(s)}^2 + \norme{0}{h(s)}^2~ds. $$

It only remains to bound $\int_0^t \norme{0}{h}^2$, and this is simple:
$$ \int_0^t \norme{0}{h(s)}^2~ds \leq C\eps^2\int_0^t \trinorme{\infty,t}{\nabla u}^2\norme{0}{R(s)}^2 + \norme{1}{\nabla u(s)}^2+\norme{1}{\rho(s)-1}^2~ds \hspace{50pt} $$
$$ \hspace{150pt} + C\int_0^t (1+\norme{\infty}{\derp{z}u_3}^2)\norme{1}{u_3}^2 + \norme{0}{F}^2~ds, $$
in which $\eps\int_0^t \norme{1}{\nabla u(s)}$ is bounded using \propref{hmee}, so we have proved the energy estimate on $R$ for $m=1$.
\newline

We now move to the $H^{m-1}_{co}$ estimates with $m\geq 2$, and once again, the main focus will be the commutators. Given $\alpha$ of length $m-1$, the equation on $Z^\alpha R$ is
$$ l\eps(\derp{t}Z^\alpha R + u\cdot\nabla Z^\alpha R) + \rho P'(\rho) Z^\alpha R = Z^\alpha h - {\cal C}^\alpha_R, $$
with ${\cal C}^\alpha_R = l\eps[Z^\alpha, u\cdot\nabla] R + [Z^\alpha, \rho P'(\rho)]R$. Repeating the above procedure, we reach the estimate
$$ c\eps\norme{0}{Z^\alpha R(t)}^2 + c\int_0^t \norme{0}{Z^\alpha R(s)}^2~ds \leq C\norme{m-1}{R(0)}^2 \hspace{100pt} $$
\begin{equation} \hspace{100pt} + C\int_0^t \eps\trinorme{\Lip,t}{U}\norme{m-1}{R(s)}^2+\norme{m-1}{h(s)}^2 + \norme{0}{{\cal C}^\alpha_R(s)}^2~ds. \label{hmndr1} \end{equation}

Controlling $\int_0^t \norme{m-1}{h}^2$ is immediate, except for the term $\rho^2 u_3 \derp{z}u_3$, which figures in $\rho^2 u\cdot\nabla u_3$. We start by replacing $\rho\derp{z}u_3$ by using (\ref{NSdiv}),
$$ \rho^2 u_3 \derp{z}u_3 = -\rho u_3(\derp{t}\rho + \rho \divt u\tang + u\tang\cdot\nabla\tang \rho + u_3\derp{z}\rho), $$
and the problematic term here is clearly $\rho u_3^2 \derp{z}\rho$. Directly using the tame estimate basically yields $\int_0^t \norme{m-1}{R}^2$, but there is no factor $\eps$ in this term to allow us to hope to absorb it. 
However, we have a factor $u_3$, and we use the assumption that, for a given $\delta>0$, there exist $\sdel{z}>0$ which does not depend on $\eps$, such that $|u_3(x)|<\delta$ for $x$ in the strip $\sdel{\omega} = \mathbb{R}^2\times [0,\sdel{z}]$. 
Now, we split the $H^{m-1}_{co}$ norm of $\rho u_3^2 R$ into two parts,
$$ \int_0^t \norme{m-1}{\rho u_3^2 R(s)}^2 ~ds = \int_0^t \norme{m-1,\sdel{\omega}}{\rho u_3^2 R(s)}^2 + \norme{m-1,\Omega\backslash \sdel{\omega}}{\rho u_3^2 R(s)}^2~ds. $$
Here, for $\omega\subset \Omega$, we set $\norme{m-1,\omega}{f(s)}$ as a sort of semi-norm of $f(s)$ restricted to $\omega$:
$$ \norme{m-1,\omega}{f(s)}^2 = \sum_{|\beta|\leq m-1} \norme{0}{(Z^\beta f)(s) |_{\omega}}^2 $$
We apply the tame estimate to both norms:
\begin{eqnarray*} \int_0^t \norme{m-1}{\rho u_3^2 R}^2 & \leq & C\int_0^t \trinorme{L^\infty(\sdel{\omega}),t}{\rho u_3^2}^2\norme{m-1}{R}^2 + \trinorme{\infty,t}{R}^2\norme{m-1}{\rho u_3^2}^2 \\
& & \hspace{50pt} + \trinorme{\infty,t}{\rho u_3^2}^2\norme{m-1,\Omega\backslash\sdel{\omega}}{\derp{z}\rho}^2~ds. \end{eqnarray*}
The two key terms are the first and the last. To deal with the last term, we use \propref{comprop} to write
$$ \norme{m-1,\Omega\backslash\sdel{\omega}}{\derp{z}\rho}^2 \leq C\sum_{|\beta|< m} \norme{0}{\derp{z} Z^\beta \rho}^2, $$
and we note that, for $z\geq \sdel{z}$, $\phi(z)\geq \phi(\sdel{z})$, therefore
$$ |\derp{z} Z^\beta \rho(x)|\leq \frac{1}{\phi(\sdel{z})}|Z_3 Z^\beta \rho(x)| $$
for $x\neq \sdel{\omega}$. This means that conormal derivatives are equivalent to standard derivatives away from the boundary, thus, $\norme{m-1,\Omega\backslash\sdel{\omega}}{\derp{z}\rho}^2 \leq \norme{H^m_{co}(\Omega)}{\rho-1}^2$. 
The first term is led by $\norme{L^\infty(\sdel{\omega})}{\rho u_3^2}^2$, which, given the boundedness of $\rho$ and the properties of $u_3$ on $\sdel{\omega}$, is bounded by $c_1^2\delta^4$, which is smaller that $\frac{c}{2}$ if $\delta$ is small enough, where $c$ is the coefficient on the left-hand side of (\ref{hmndr1}): we therefore choose $\delta$ so that this term can be absorbed.

The other terms in $h$ are straight-forward, as we can use \estiref{hmee} on the $\eps^2\int_0^t \norme{m}{\nabla u}^2$ term that comes from the order-two terms of $h$ - and the factor $\eps^2$ is essential, as it leaves a factor $\eps$ which will allow us to close the complete estimate for $\eps$ and $t$ small. In total, we get
$$ \int_0^t \norme{m-1}{h(s)}^2~ds \leq C\norme{m}{U(0)}^2 + \left(\eps Q(1+\trinorme{\Lip,t}{U}^2) + \frac{c}{2}\right)\int_0^t \norme{m-1}{R(s)}^2~ds $$
\begin{equation} \hspace{50pt} + C\int_0^t Q(1+\trinorme{\Lip,t}{U}^2 + \trinorme{1,\infty,t}{\nabla u}^2) (\norme{m}{U(s)}^2+\norme{m-1}{\nabla u\tang(s)}^2+\norme{m-1}{F(s)}^2)~ds \label{hmndr2} \end{equation}

The estimation of the commutators is also mostly straight-forward. $[Z^\alpha,u\cdot\nabla]R$ can be controlled by using \propref{comlem}, bearing in mind that there is a factor $\eps$ in front of it:
$$ \eps^2\int_0^t \norme{0}{[Z^\alpha,u\cdot\nabla]R}^2~ds \leq C\eps \int_{0}^{T}(\trinorme{\Lip,t}{u}^{2}+\trinorme{1,\infty,t}{\nabla u}^{2})\norme{m-1}{R(s)}^{2}~ds \hspace{30pt} $$
$$ \hspace{110pt} +C\eps \int_{0}^{t}\trinorme{1,\infty,t}{R}^{2}\norme{m-1}{\nabla u(s)}^{2}+\trinorme{\infty,t}{R}^{2}\norme{m}{u(s)}^{2}~ds. $$
The other commutator, $[Z^\alpha,\rho P'(\rho)]$ does not have a factor $\eps$, so we need to gain a derivative on $R$ by using \propref{funcprop},
$$ \int_{0}^{t}\norme{0}{[Z^{\alpha},\rho^\gamma]R}^{2}~ds \leq C\trinorme{1,\infty,t}{\rho}^{2\gamma}\int_{0}^{t}\norme{m-2}{R(s)}^{2}~ds+\trinorme{\infty,t}{R}^{2}\int_0^t \norme{m-1}{(\rho-1)(s)}^{2}~ds , $$
and use (\ref{l2time}) to extract a factor $t$ in the first term. This leaves us with an isolated $t\trinorme{1,\infty,t}{\rho}$, to which we apply the anistropic Sobolev embedding, \propref{sobo}. Thus,
$$ t\trinorme{1,\infty,t}{\rho}^2 \leq Ct\left(1+\norme{4}{R(0)}^2 +\trinorme{4,t}{U}^2+ \int_0^t \norme{5}{R(s)}^2~ds\right) $$
Note that it is this last inequality that restricts us to $m-1\geq 5$. We conclude the proof of \estiref{hmeeR} by combining (\ref{hmndr1}), (\ref{hmndr2}) and these bounds on the commutators. $\square$

\subsection{$L^\infty$ estimates}

As stated in the introduction, we will control $L^\infty$ norms of $R$ with $L^2$-in-time bounds by virtue of (\ref{l2timeinf}),
\begin{equation} \trinorme{1,\infty,t}{R}^2 \leq \norme{1,\infty}{R(0)}^2 + C\int_0^t \norme{1,\infty}{\derp{t}R(s)}^2+\norme{1,\infty}{R(s)}^2~ds . \label{l2timeinfr} \end{equation}
Let us define $Y_t$ the space of functions $f$ satisfying $\norme{Y_t}{f}^2 := \int_0^t \norme{1,\infty}{f(s)}^2+\norme{1,\infty}{\derp{t}f(s)}^2~ds$ finite.

\begin{esti} Under the conditions of \assuref{assu1}, there exists an increasing, positive function $Q$ on $\rplus$ such that, for $t\leq T^*$ and $\eps\leq \eps_0$,
$$ \norme{Y_t}{R}^2 \leq Q(M_0) + (t+\eps)Q(M+M_F) $$ \label{infR} \end{esti}
This ends the proof of \thref{thEE} providing we can pick up \assuref{assu1}.
\newline

\preu{Proof:} the main tool in this proof will be the Duhamel formula for the ordinary differential equation $\eps f' + \til{p}f = \til{h}$. We reach this ODE by considering $R$ along the characteristics of the transport equation $\derp{t}R+u\cdot\nabla R = 0$, in other words
$$ f(t,x) = R^X(t,x) := R(t,X(t,x)), $$
where $X(t,x)$ satisfies $\derp{t}X(t,x) = u(t,X(t,x))$ and $X(0,x)=x$. We extend the notation $g^X$ to any function $g$ as above. 
Thus, we have the identity $\derp{t}(R^X) = (\derp{t}R+u\cdot\nabla R)^X$, so (\ref{eqnR}) becomes,
$$ \eps l(\rho^X) \eps\derp{t}(R^X) + \rho^X P'(\rho^X) R^X = h^X. $$

For the higher-order estimates, it is important to apply the conormal derivatives first, then follow the flow of $u$. So, the equation we are interested in is
$$ \eps l(\rho^X)\eps\derp{t}((Z^\alpha R)^X) + \rho^X P'(\rho^X) (Z^\alpha R)^X = (Z^\alpha h)^X - ({\cal C}_{R}^\alpha)^X, $$
with ${\cal C}_R^\alpha=l\eps[Z^\alpha,u\cdot\nabla]R+ [Z^\alpha,\rho P'(\rho)]R$ as in the previous paragraph, and $\alpha$ is either of length $\leq 1$, or of length 2 with $\alpha_0\geq 1$ (these are the $\alpha$ that intervene in $Y_t$).

To lighten the load, we introduce the following notations: $g_\alpha = Z^\alpha h - {\cal C}^\alpha_R$,
$$ j(s',s,x) = \int_s^{s'} \frac{[\rho^X P'(\rho^X)](\sigma,x)}{\eps l(\rho^X)(\sigma,x)} ~d\sigma ~,~\mathrm{and}~ J(s,x) = \int_0^s \frac{g_\alpha^X (\sigma,x)}{\eps l(\rho^X)(\sigma,x)} e^{-j(s,\sigma,x)}~d\sigma. $$
The Duhamel formula for this equation then reads:
$$ (Z^\alpha R)^X(s,x)=Z^\alpha R(0,x)e^{-j(s,0,x)} + J(s,x) . $$
We integrate the square of this equality in time between $0$ and $t$, which yields
\begin{equation} \int_0^t [(Z^\alpha R)^X(s,x)]^2~ds \leq C\left[\norme{\infty}{Z^\alpha R(0)}^2\int_0^t e^{-2j(s,0,x)}~ds + \int_0^t J^{2}(s,x)~ds\right] . \label{duhamel2} \end{equation}
Again, we use the uniform bounds on $\rho$ to get, for any $s,~s' \in \rplus$, $x\in\Omega$,
\begin{equation} j(s',s,x) \geq \int_s^{s'} \frac{c}{\eps} = \frac{c(s'-s)}{\eps}. \label{jbelow} \end{equation}
Thus,
\begin{equation} \int_0^t e^{-cj(s,0,x)}~ds \leq C\eps \label{jest} , \end{equation}
which deals with the first term of the right-hand side of (\ref{duhamel2}). For the second term, we use (\ref{jbelow}) again and write the integral to reveal a convolution in the time variable:
\begin{eqnarray*} \int_0^t J^2(s,x)~ds & \leq & \int_0^t \left[\int_\mathbb{R} |g_\alpha^X|(\sigma,x)\mathbf{1}_{(0,t)}(\sigma) \frac{e^{-c\eps^{-1}(s-\sigma)}}{\eps}\mathbf{1}_{(0,t)}(s-\sigma)~d\sigma\right]^2 ~ds \\
 & \leq & \norme{L^2(0,t)}{\left(g_\alpha^X(\cdot,x)\mathbf{1}_{(0,t)}\right) \ast \left(\eps^{-1}e^{-c\eps^{-1}(\cdot)}\mathbf{1}_{(0,t)}\right)}^2 \\
 & \leq & \norme{L^2(0,t)}{g_\alpha^X(\cdot,x)}^2 \norme{L^1(0,t)}{\eps^{-1}e^{-c\eps^{-1}(\cdot)}}^2 \\
 & \leq & C\int_0^t \norme{\infty}{g_\alpha (s)}^2~ds, \end{eqnarray*}
by (\ref{jest}) and standard convolution inequalities.
\newline

So now we only need to control the $Y_t$ norm of $h$ and the $L^\infty$ norm of the commutators. We remind the reader that
\begin{eqnarray*} h & = & \eps[\rho(\mu\Delta\tang u_3+\lambda\derp{z}\divt u\tang+\sigma(\nabla U)_3) \\
 & & \hspace{20pt} -l\derp{z}u\tang\cdot\nabla\tang\rho-l(\div u+\derp{z}u_3)R]-\rho^2(\derp{t}u_3+u\cdot\nabla u_3)+\rho F_3. \end{eqnarray*}
The starting point here is to notice that, for two functions $f$ and $g$,
\begin{equation} \norme{2,\infty}{fg} \leq C(\norme{2,\infty}{f}\norme{\infty}{g} + \norme{1,\infty}{f}\norme{2,\infty}{g}). \label{fg2inf} \end{equation}
In the case of $h_1:=\rho(\mu\Delta\tang u_3 + \lambda\derp{z}\divt u\tang)$, $f=\rho$, so
$$ \norme{2,\infty}{h_1}^2 \leq C(\norme{2,\infty}{\rho}^2(\norme{2,\infty}{u}^2+\norme{1,\infty}{\nabla u\tang}^2) + \norme{1,\infty}{\rho}(\norme{4,\infty}{u}^2 + \norme{3,\infty}{\nabla u\tang}^2)), $$
in which we apply the Sobolev inequality to all terms except $\norme{1,\infty}{\nabla u\tang}^2$ (it is part of the total quantity we wish to bound), and use the Young inequality to split the products up. 
In the process, we obtain $(\eps^2 \int_0^t \norme{4}{\derp{z}\nabla u\tang(s)}^2~ds)^2$, which is controlled by using \estiref{hmvort}. Therefore, in total, using all the estimates available to us,
$$ \int_0^t \eps^2\norme{2,\infty}{h_1}^2 ~ds \leq \eps Q(M_0) + \eps \int_0^t Q(1+\norme{7}{U}^2+\norme{6}{\nabla u\tang}^2+\norme{6}{R}^2 +\norme{1,\infty}{\nabla u\tang}^2 + \norme{1,\infty}{R}^2)~ds, $$
\begin{equation} \mathrm{thus}~~\int_0^t \eps^2\norme{2,\infty}{h_1}^2 ~ds \leq Q(M_0) + C\eps(t+1)Q(M+M_F) + C\eps t \trinorme{1,\infty,t}{R}^2. \label{densinfh1} \end{equation}
Likewise, setting $h_2:= \derp{z}u\tang\cdot\nabla\tang\rho+(\div u + \derp{z}u_3)R + \sigma_3 = \derp{z}u\cdot\nabla \rho + R\div u + \sigma_3$, we use (\ref{fg2inf}) again, applied to the space $Y_t$ instead of $W^{2,\infty}_{co}$, and we get, after applying the Sobolev embedding to all the terms except $\norme{1,\infty}{\nabla u\tang}^2$, and $\norme{Y_t}{R}^2$, to get
$$ \eps^2\norme{Y_t}{h_2}^2 \leq Q(M_0) + \eps(t+1)Q(M+M_F) + C\eps^2 \norme{Y_t}{R}^2 + C\eps^2 \int_0^t \norme{6}{\derp{zz}u_3}^2~ds . $$
We have once again used \estiref{hmvort} to control the $H^4_{co}([0,T]\times\Omega)$ norm of $\eps \derp{z}\nabla u\tang$, and we have used (\ref{l2time}) to get an $L^2$-in-time norm of $\derp{z}u_3$ (which can then be controlled using \propref{hmee}) and $\derp{zz}u_3$. 
Thus, we replace $\eps\derp{zz}u_3$ by its expression in (\ref{NSrwR}), use the tame estimate, take the supremum inside the integral and extract a factor $t$, and this yields
\begin{equation} \eps^2\norme{Y_t}{h_2}^2 \leq Q(M_0) + \eps(t+1)Q(M+M_F) + C\eps \norme{Y_t}{R}^2. \label{densinfh2} \end{equation}

In the source term $h$, it remains to examine $h_3:= \rho^2(\derp{t}u_3+u\cdot\nabla u_3+F_3)$. The difference here is that there is no factor $\eps$ ready to provide us with a small parameter. Instead, we once again extract $t$ from the integral to use as the small parameter. 
By simply using the Sobolev embedding inequality and breaking the result down almost completely with the tame estimate, we get
\begin{eqnarray} \int_0^t \norme{2,\infty}{h_3(s)}^2~ds & \leq & Q(M_0)+ C\int_0^t \norme{5}{h_3(s)}^2+\norme{4}{\derp{z}h_3(s)}^2~ds \nonumber \\
 & \leq & Q(M_0) + \int_0^t Q(1+\trinorme{1,\infty,t}{U}^2+\trinorme{1,\infty,t}{\derp{z}U}^2+M_F) \nonumber \\
 & & \hspace{35pt} \times \left(1+\norme{6}{U}^2+\norme{5}{\derp{z}U}^2+\norme{4}{u_3\derp{zz}u_3}^2~ds\right) . \label{h31} \end{eqnarray}
The only product we do not split using the tame estimates and the Young inequality is $u_3\derp{zz}u_3$, which appears in $\derp{z}h_3$. Indeed, we need the factor $u_3$ to be able to compensate for the two $z$-derivatives. We have already used in (\ref{h31}) the fact that
$$ \trinorme{\infty,t}{u_3\derp{zz}u_3}^2 \leq \trinorme{1,\infty,t}{\derp{z}u_3}^4 , $$
due to (\ref{u3phi}), and $\norme{3}{u_3\derp{zz}u_3}$ is dealt with in same way as in the proof of (\ref{comfinal}) in \propref{comlem}: multiply and divide by $\phi$, which means we are actually looking at $\norme{3}{\phi^{-1}u_3 Z_3\derp{z}u_3}$, and use the Hardy inequality to get that this quantity satisfies (\ref{comfinal}), with $g=u_3$ and $f=\derp{z}u_3$, so
$$ \int_0^t \norme{2,\infty}{h_3(s)}^2 ~ds \leq Q(M_0) + Q(1+\trinorme{1,\infty,t}{U}^2+\trinorme{1,\infty,t}{\derp{z}U}^2)\int_0^t \norme{6}{U}^2+\norme{5}{\derp{z}U}^2~ds . $$
Now we extract $t$ and use (\ref{l2time}) to get $L^2$-in-time norms on $R$ and $\derp{z}u_3$, so
\begin{equation} \int_0^t \norme{2,\infty}{h_3(s)}^2~ds \leq Q(M_0) + tQ(M+M_F) . \label{densinfh3} \end{equation}

Finally, we examine the $L^{\infty}$ norms of the commutators, ${\cal C}^{\alpha}_R$, with $|\alpha|= 1$, or $|\alpha|=2$ and $\alpha_0\geq 1$. We begin with the case $|\alpha|=1$, so $Z^\alpha=Z_j$ for a certain $j$, and
$$ {\cal C}^{\alpha}_{R} = l\eps[Z_j,u\cdot\nabla] R + [Z_j,\rho P'(\rho)]R = l\eps(Z_j u)\cdot\nabla R + l\eps u_3[Z_j,\derp{z}]R + \gamma(Z_j P)R , $$
in which the second term is either 0 ($j\neq 3$) or $-l\eps\phi' u_3\derp{z}R$. Bounding this is straight-forward, using (\ref{hardyinf}) along the way:
\begin{equation} \int_0^t \norme{\infty}{{\cal C}^\alpha_R(s)}^2~ds \leq Ct(1+\eps^2)\trinorme{1,\infty,t}{\nabla u}^2\trinorme{1,\infty,t}{R}^2 . \label{densinfc1} \end{equation}
The second case, $|\alpha|=2$ and $\alpha_0>0$ is also simple. We can write $Z^\alpha = Z_j \derp{t}$ for a certain $j$, so
\begin{eqnarray*} {\cal C}^{\alpha}_R & = & l\eps\left((Z_j\derp{t} u)\cdot\nabla R + (Z_ju)\cdot\nabla(\derp{t}R) + (\derp{t}u)\cdot\nabla (Z_j R) + (\derp{t}u_3)[Z_j,\derp{z}]R \right. \\
 & & \hspace{25pt} \left. + u_3[Z_j,\derp{z}]\derp{t}R\right) + \gamma((Z_j\derp{t}P)R+(Z_j P)(\derp{t}R)+(\derp{t} P)(Z_jR)) \end{eqnarray*}
The $L^\infty$ norm of this is, for most terms, bounded using only the Sobolev embedding and/or (\ref{hardyinf}), wherever $u_3\derp{z}$ or conormal derivatives of $u_3 \derp{z}$ appear. One can then take the supremum in time inside the integral and integrate to get a factor $t$. For instance, the last term satisfies
\begin{eqnarray*} \norme{\infty}{(\derp{t}r)(Z_jR)}^2 & \leq & C(\norme{2,\infty}{u\tang}^2+\norme{1,\infty}{\derp{z}u_3}^2+\norme{1,\infty}{\rho}^2)\norme{1,\infty}{R}^2 \\
 & \leq & C(1+\norme{5}{U}^2+\norme{4}{\nabla u\tang}^2+\norme{3}{R}^2+\norme{1,\infty}{\derp{z}u_3}^2)\norme{1,\infty}{R}^2 , \end{eqnarray*}
then we apply \estiref{hmnd} (a) on $\norme{1,\infty}{\derp{z}u_3}$. Two specificities do appear though in the terms containing $u_3 \derp{z}R$ and the like.
\begin{itemize}
\item The term $l \eps u_3[Z_3,\derp{z}]\derp{t}R = -l\eps\phi' u_3\derp{zt}R$ leads to two conormal derivatives on $R$. We cannot extract a factor $t$ from the integral here, since we want $L^2$-in-time norms of $\norme{1,\infty}{\derp{t}R}$. However, we do have a factor $l\eps$, which will act as the small parameter. So,
$$ \int_0^t \norme{\infty}{l\eps\phi' u_3\derp{zt}R}^2~ds \leq C\eps^2\trinorme{\infty,t}{\derp{z}u_3}^2\int_0^t \norme{1,\infty}{\derp{t}R(s)}^2~ds . $$
\item The term $l \eps(Z_j\derp{t}u)\cdot\nabla R$ contains $\eps(Z_j\derp{t}u_3)\derp{z}R$, which, after using (\ref{hardyinf}), leads to a term with two conormal derivatives on $\eps\derp{z}u_3$ (and no more than one on $R$, which is therefore not problematic), so we re-use the trick we used when estimating $\eps h_2$: we use the Sobolev inequality on the term with two conormal derivatives, to get
$$ \int_0^t \eps^2\norme{2,\infty}{\derp{z}u_3(s)}^2~ds \leq M_0 + C\int_0^t \eps^2\norme{5}{\derp{z}u_3(s)}^2 + \norme{5}{\eps\derp{zz}u_3(s)}^2~ds , $$
and we replace $\eps\derp{zz}u_3$ by its expression in (\ref{NS}), use the tame estimate, take the supremum inside the integral and extract a factor $t$.
\end{itemize}

In total, combining the above consideration with inequalities (\ref{densinfh1}) to (\ref{densinfc1}), we conclude that
$$ \int_0^t \norme{\infty}{g_\alpha(s)}^2~ds \leq Q(M_0) + (t+\eps) Q(M+M_F) , $$
which ends the proof of the estimate. $\square$

\section[Proof of \thref{thEE}, part V: conclusion]{Proof of \thref{thEE}, part V \\ Conclusion
 \sectionmark{Proof of \thref{thEE} part V: conclusion}}
\sectionmark{Proof of \thref{thEE} part V: conclusion}

Let us consider times $t\leq T'$, where
$$ T' = \sup\{ t\leq T_0 ~|~ {\cal E}_m(t,U) + \trinorme{1,\infty,t}{\derp{z}u_3}^2 \leq Q(2M_0)\} . $$
The time $T'$ depends \textit{a priori} on $\eps$. With this, we pick up the uniform boundedness from below of $\rho$ and the uniform smallness of $u_3$ near the boundary that we have so far assumed.
\begin{itemize}
\item Let $\rho(0,x)\geq c_{0}$. Equation (\ref{NSdiv}) provides us with the differential equation
$$ \derp{t}\rho^{X}(t,x) = - (\rho\div u)^{X}(t,x), $$
where $f^{X}$ is once again $f$ following the flow of $u$. Thus,
$$ \rho^{X}(t,x) = \rho(0,x)\exp\left(-\int_{0}^{t}\div u(s,x)~ds\right), $$
and $\div u$ is uniformly bounded on $[0,T']\times \Omega$, so
$$ |\rho(t,x)|\geq c_{0}e^{-Q(2M_0)T'} := c_0', $$
for $t\leq T'$, and $c_0'>0$ can be chosen independent of $\eps$ (if $T'=+\infty$ for some $\eps$, we replace with some finite $T'_0$ independent of $\eps$). The density $\rho$ is therefore uniformly bounded from below on $[0,T']$.
\item For $t\leq T'$, we can repeat the proof of (\ref{linfnd1}),
$$ \norme{\infty}{\derp{z}u_3}^2 \leq Q(1+\norme{\Lip}{U(0)}^2) + tQ(1+\trinorme{4,t}{U}^2+\trinorme{3,t}{\derp{z}U}^2 + \trinorme{1,\infty,t}{\derp{z}\rho}^2), $$
which is only a reading of (\ref{NSdiv}) combined with the Sobolev embedding and property (\ref{l2timeinf}), hence $\derp{z}u_3(t)$ has a uniform bound for $t\leq T'$, thus $|u_3(t,z)|\leq M'z$ for some $M'$ when $t\leq T'$. 
As a result, we get uniform smallness of $u_3$ near the boundary: $|u_3(t,z)|\leq M'\delta$ when $z<\delta$. \end{itemize}

For $t\leq T'$, we have the bounds required to make the \textit{a priori} estimation process valid, so we have
$$ {\cal E}_m(t,U) \leq Q(M_0)+(t+\eps)Q(Q(2M_0)+M_F). $$
For $\eps\leq \eps_0$ small enough, we see that the right-hand side is smaller than $Q(2M_0)$ for $t\leq T^*\leq T'$, with $T^*$ depending on $\eps_0$ but not on $\eps$. Thus $T^*$ can be chosen independent of $\eps$, and the theorem is proved. $\square$
\vspace{12pt}

\textbf{\underline{Acknowledgements.}} This is one of the results from my Ph.D thesis, prepared at the University of Rennes 1. Therefore, I am deeply grateful to my supervisor Fr\'ed\'eric Rousset for the opportunity to work on this topic, for his guidance and patience throughout the preparation of the present article. 
I would also like to thank Mark Williams for the discussions that have helped improve the clarity of the paper.

\begin{small}
\bibliography{nsebib}
\bibliographystyle{abbrv}
\end{small}

\end{document}